\documentclass{article}

\usepackage{float}
\usepackage{PRIMEarxiv}
\usepackage[T1]{fontenc}    
\usepackage{hyperref}       
\usepackage{url}            
\usepackage{booktabs}       
\usepackage{amsfonts}       
\usepackage{nicefrac}       
\usepackage{microtype}      
\usepackage{lipsum}
\usepackage{fancyhdr}       
\usepackage{graphicx}       
\usepackage{amsmath}
\graphicspath{{media/}}     
\usepackage{subcaption}
\numberwithin{equation}{section}
\usepackage{amsthm}

\newtheorem{proposition}{Proposition}[section]
\usepackage{enumitem}
\title{
Task-Contribution Dynamics of Adam in Physics-Informed Neural Networks: Slow-Manifold Analysis and Gradient Modulation
}

\author{
  Enzhan Zheng \\
  School of Mathematical Sciences \\
  Dalian University of Technology  \\
  Dalian, Liaoning 116024, China \\
  zhengenzhan@mail.dlut.edu.cn \\
  \And
  Xiaoliang Song\thanks{
  	Corresponding author:
  	\texttt{songxiaoliang@dlut.edu.cn}
  } \\
  School of Mathematical Sciences \\
  Dalian University of Technology \\
  Dalian, Liaoning 116024, China \\
  \texttt{songxiaoliang@dlut.edu.cn} \\
}

\begin{document}
\maketitle

\begin{abstract}
	Physics-informed neural networks (PINNs) often suffer from imbalanced optimization of competing loss components, particularly when approximating solutions with high-frequency or multiscale features. Although existing balancing strategies have proven effective in alleviating this difficulty, the mechanisms governing the evolution of individual task contributions to Adam's parameter updates remain insufficiently understood. We introduce task-wise indicators to characterize this evolution and identify three distinct training regimes: smooth variation, intermittent spikes, and persistent oscillations. Focusing on the early smooth regime, we combine discrete-time analysis with a continuous-time singular perturbation framework under the full-batch setting and the parameter regime $(1-\beta_2)/(1-\beta_1)\ll1$. A first-order approximation of the locally attracting slow manifold reveals two factors governing the step-to-step redistribution of task contributions: curvature-related alignment disparities and differences in the coordinatewise distributions of preconditioned task gradients. Both factors involve Adam's shared second-moment estimate. Inspired by AdaTask, we propose Gradient-Modulated Adam (GMAdam), which modulates incoming task gradients using task-specific second-moment estimates. Numerical experiments on high-frequency Poisson and wave equations, the Allen--Cahn equation, and two-dimensional lid-driven cavity flow ($Re=100,400$) demonstrate consistent improvements in predictive accuracy over AdaTask at comparable computational cost.
\end{abstract}

\keywords{multi-objective optimization \and Adam optimizer \and physics-informed neural networks \and deep learning}

\section{Introduction}
The universal approximation capabilities of deep neural networks (DNNs) have given rise to a new paradigm for scientific computing involving partial differential equations (PDEs)  \cite{E_2017,han2018solving,e2017deepritzmethoddeep,liao2019deep,zang2020weak}. In this context, physics-informed neural networks (PINNs) \cite{raissi2019physics} have emerged as a representative methodology, in which observational data and governing physical laws are jointly incorporated into the training process. Over the past several years, PINNs have facilitated a series of promising results across various fields including fluid mechanics \cite{sun2020surrogate,cai2021physics,eivazi2022physics}, geosciences \cite{haghighat2021sciann,smith2022hyposvi}, bio-engineering \cite{kissas2020machine,sahli2020physics}, electromagnetics \cite{9403414,kovacs2022conditional}, and meta-material design \cite{liu2019multi,chen2020physics}.

However, despite the notable empirical success of PINNs, significant challenges persist when they are built upon conventional neural network architectures. One fundamental challenge arises from the competing objectives in PINNs training, whose gradients may exhibit substantial disparities in magnitude and conflicts in direction \cite{wang2025gradient}. These issues can be further exacerbated in the presence of high-frequency and multi-scale features of the underlying PDE solution, thereby impeding convergence toward parameters that yield low values across all loss components. To mitigate gradient magnitude imbalance, \cite{wang2021understanding} proposed a periodically updated reweighting strategy (Grad Norm) aimed at equalizing the gradient norms of individual loss terms. A complementary perspective was provided by Wang et al.~\cite{wang2022and}, who investigated the convergence behavior of individual PINN loss components within the neural tangent kernel (NTK) framework \cite{jacot2018neural}. Building on this analysis, they developed an NTK-based reweighting strategy to promote balanced decreases across the different loss components.

While the preceding reweighting strategies are motivated by the behavior of task gradients, another approach seeks to address gradient magnitude imbalance through Adam itself, one of the most widely used optimizers in PINN training. In the context of multi-task learning with Adam, Yang et al.~\cite{AdaTask_AAAI2023} introduced the relative accumulated update (rAU) metric to characterize task dominance at the level of individual parameter updates, linking this phenomenon to imbalances in accumulated gradient statistics. This observation motivated the development of AdaTask, which employs task-specific adaptive learning rates to mitigate such dominance. 
For PINNs, Yao et al.~\cite{pmlr-v202-yao23c} presented MultiAdam, which likewise employs task-specific adaptive learning rates but averages the resulting updates rather than summing them as in AdaTask.

Despite these developments, the mechanisms governing the evolution of task contributions under standard Adam remain insufficiently understood. Although NTK-based analyses have provided insights into the convergence behavior of individual PINN loss components under gradient descent (GD) \cite{wang2022and}, Adam exhibits distinct dynamics arising from momentum accumulation and coordinatewise preconditioning. From a dynamical-systems perspective, Barakat and Bianchi \cite{barakat2021convergence} established a continuous-time ODE approximation for stochastic Adam, whereas Ma et al.~\cite{pmlr-v145-ma22a} examined the qualitative dynamics of full-batch Adam, including oscillations and spikes in the training loss. Nevertheless, these theoretical advances do not directly explain how Adam's shared second-moment estimate influences the relative contributions of individual loss components to parameter updates. In particular, the mechanisms governing the step-to-step redistribution of these contributions remain to be characterized.

In this work, we investigate the redistribution of task contributions under Adam, focusing on the early smooth phase of training. Our theoretical analysis combines discrete-time analysis with a continuous-time dynamical formulation in the full-batch setting and under the parameter regime $\varepsilon=(1-\beta_2)/(1-\beta_1)\ll1$, revealing the factors governing the step-to-step evolution of relative task contributions. Motivated by these findings, we further develop a task-aware gradient-modulation scheme inspired by AdaTask \cite{AdaTask_AAAI2023}. The main contributions of this work are summarized as follows:

\begin{enumerate}[label=(\roman*), leftmargin=*]
	
	\item We introduce task-wise indicators to quantify the relative contributions of individual loss components to Adam's parameter updates and identify three characteristic regimes: smooth evolution, intermittent spikes, and persistent oscillations, motivating our focus on the early smooth regime.
	
	\item Under the full-batch setting and $(1-\beta_2)/(1-\beta_1)\ll1$, we establish a continuous-time correspondence for Adam and derive a locally attracting slow manifold using geometric singular perturbation theory. The resulting first-order approximation characterizes task-contribution redistribution through curvature-related alignment disparities and differences in the coordinatewise distributions of preconditioned task gradients.
	
	\item We examine how task dominance in the shared second-moment estimate can limit the redistribution of update contributions. Motivated by an equivalent formulation of AdaTask, we propose Gradient-Modulated Adam (GMAdam), which modulates task gradients using task-specific second-moment estimates.
	
	\item We evaluate the proposed method on four PDE benchmarks: the Poisson, wave, and Allen--Cahn equations, and lid-driven cavity flow ($Re=100,400$). Comparisons with AdaTask and two loss-reweighting methods demonstrate consistently higher predictive accuracy than AdaTask at comparable computational cost.
	
\end{enumerate}

The remainder of this paper is organized as follows. Section~2 reviews the PINN formulation and representative strategies for mitigating training imbalance. Section~3 develops the discrete- and continuous-time analyses of Adam's task-contribution dynamics, followed by the construction of the proposed optimization scheme. Section~4 presents numerical experiments on the four benchmark problems and compares the proposed method with existing approaches.

\section{Preliminaries}
\label{sec:Preliminaries}
\subsection{Physics-informed Neural Networks}
\label{subsec:PINNs}
In this section, we begin with a concise review of physics-informed neural networks (PINNs) \cite{raissi2019physics}. Let us consider the following partial differential equation (PDE) problem posed on a bounded domain: $\Omega \subset \mathbb{R}^{p_0}$
\begin{align}
	&\mathcal{N}[u]=f(x), \quad x \in \Omega \label{eq:pde}\\
	&\mathcal{B}[u]=g(x), \quad  x \in \partial \Omega \label{eq:bc} 
\end{align}
where $\mathcal{N} [ \cdot ] $ is a differential operator and $\mathcal{B}[\cdot]$ represents the associated boundary operator. Additionally, $u:\bar{\Omega} \rightarrow \mathbb{R}$ refers to the unknown solution of the governing PDE system \eqref{eq:pde}--\eqref{eq:bc}. In the time-dependent setting, the variable $x=(\mathbf{x},t)$ denotes the space--time coordinate and, accordingly, the initial conditions are subsumed into the boundary operator $\mathcal{B}[\cdot]$. 

PINNs are designed to approximate the desired solution $u(x)$ by a neural network function $u_\theta(x)$, with the tunable parameters $\theta$ defining the underlying representation. In the original PINNs formulation \cite{raissi2019physics}, $u_\theta(x)$ is instantiated as multi-layer perceptrons (MLPs), a choice that has remained prevalent across subsequent PINNs studies. For an $L$-layer MLP, we consider the corresponding function $u_{\theta}:\mathbb{R}^{p_0} \rightarrow \mathbb{R}$, whose input coordinates $x$ are first propagated through $L-1$ layers according to
\begin{equation*}
	\begin{aligned}
		&\bar{h}_l(x) = W_l h_{l-1}(x)+b_l\\
		&h_{l}(x) = \phi(\bar{h}_l(x))\\
		&h_{0}(x) = x
	\end{aligned}
\end{equation*}
for $l=1,\cdots,L-1$. Here $W_l \in \mathbb{R}^{p_l \times p_{l-1}}$ and $b_l \in \mathbb{R}^{p_l}$ denote the weight matrix and bias vector of the $l$-th hidden layer; $\phi:\mathbb{R}\rightarrow \mathbb{R}$ is an activation function applied componentwise. Then the final network output takes the form:
\begin{equation*}
	u_{\theta}(x)=\bar{h}_L(x) = W_L \cdot h_{L-1}(x)+b_L
\end{equation*}
where $W_L\in\mathbb{R}^{1\times p_{L-1}}$ and $b_L\in\mathbb{R}^{1}$. Consequently, the trainable parameter set can be written as $\theta := \{ W_l, b_l \}_{l=1}^L$ with the total number of trainable degrees of freedom equal to $d$. The weight matrices $\{W_l\}_{l=1}^L$ are initialized via the Glorot scheme \cite{glorot2010understanding}, and the bias vectors $\{b_l\}_{l=1}^L$ are sampled from a standard normal distribution.

The training of $\theta$ can be cast as the following multi-objective optimization problem.
\begin{equation}
	\label{pro:PINNs}
	\min \mathcal{L}(\theta) = \mathcal{L}_{bc}(\theta) + \mathcal{L}_{r}(\theta)
\end{equation}
where
\begin{align*}
	&\mathcal{L}_{bc}(\theta)=\frac{1}{N_{bc}}\sum\limits_{i=1}^{N_{bc}} \vert \mathcal{B}[u_{\theta}](x^i_{bc})-g(x^i_{bc}) \vert^2,\\
	&\mathcal{L}_{r}(\theta)=\frac{1}{N_r}\sum\limits_{i=1}^{N_r} \vert \mathcal{N}[u_{\theta}](x^i_r)-f(x^i_r) \vert^2.
\end{align*}
Here, $\mathbf{x}_{bc}=\{x_{bc}^i\}_{i=1}^{N_{bc}}$ and $\mathbf{x}_r=\{x_r^i\}_{i=1}^{N_r}$ denote the boundary and interior collocation sets, comprising $N_{bc}$ and $N_r$ points sampled from their respective domains. The resulting empirical objective is then typically minimized using a gradient-based optimizer such as Adam \cite{kingma2015iclr-adam}. By construction, the PDE residual loss $L_{r}$ measures how severely the neural network prediction violates the considered equation \eqref{eq:pde} at interior collocation points. Only by simultaneously minimizing the boundary and initial condition losses as the PDE residual loss decreases can the network be guided toward the solution of equation \eqref{eq:pde}. This emphasizes the necessity of coordinated optimization across all loss components. Unfortunately, in practice, the training dynamics of PINNs often fail to exhibit such balanced reductions, particularly when the target solution contains high-frequency and multi-scale features \cite{wang2021understanding,wang2022and}.

\subsection{Strategies for Mitigating Gradient Magnitude Imbalance}
This observation motivates a closer examination of the optimization mechanisms underlying such unbalanced training dynamics. Previous studies \cite{wang2021understanding,wang2022and} have identified pronounced disparities in the gradient magnitudes of different loss components as a characteristic feature of practical PINN training. Such gradient imbalances may cause a loss term with disproportionately large gradients to dominate the overall descent direction, thereby biasing the optimization trajectory toward satisfying only a subset of the imposed constraints. A natural remedy is to assign appropriate weights to individual loss terms to regulate their contributions throughout training. Following this idea, the original PINN objective \eqref{pro:PINNs} can be rewritten as a weighted composite loss:
\begin{equation}
	\min \mathcal{L}(\theta)
	=
	\lambda_{bc}\mathcal{L}_{bc}(\theta)
	+
	\lambda_{r}\mathcal{L}_{r}(\theta).
\end{equation}
However, even for a fixed PDE problem, appropriate weighting coefficients may vary across training stages owing to the evolving optimization landscape. This motivates adaptive weighting mechanisms, among which gradient-norm-based weighting (Grad Norm) \cite{wang2021understanding} and neural tangent kernel (NTK)-based weighting \cite{wang2022and} are two representative strategies widely adopted in PINN training.

For the subsequent discussion, we adopt the implementations of these two weighting strategies described by Wang et al.~\cite{wang2023expertsguidetrainingphysicsinformed}. While the two methods rely on distinct diagnostic quantities, they share a common smoothing rule for updating the loss weights. Specifically, the candidate weights $\hat{\lambda}_{new}$, newly estimated via the corresponding imbalance indicators, are incorporated through an exponential moving average update to dampen fluctuations arising from stochastic sampling.
\begin{equation*}
	\label{eq:alpha}
	\lambda_{new} = \alpha \lambda_{old} + (1-\alpha) \hat{\lambda}_{new}
\end{equation*}

In the case of Grad Norm, the loss weights are adjusted to bring the norms of the weighted gradients toward a common target value. Since successful training is expected to drive the gradient norms associated with the individual loss terms asymptotically toward zero, the reference level must track their time-varying amplitude. This consideration suggests the following update criterion:
\begin{equation*}
	\begin{aligned}
		\Vert \hat{\lambda}_{bc}
		\nabla_{\theta}\mathcal{L}_{bc}(\theta_t)\Vert=
		\Vert \hat{\lambda}_{r}
		\nabla_{\theta}\mathcal{L}_{r}(\theta_t)\Vert=
		\Vert\nabla_{\theta}\mathcal{L}_{bc}(\theta_t)\Vert
		+
		\Vert\nabla_{\theta}\mathcal{L}_{r}(\theta_t)\Vert.
	\end{aligned}
\end{equation*}

In contrast, the NTK-based weighting scheme draws upon NTK theory \cite{jacot2018neural}. In the gradient-flow limit, the evolution of pointwise constraint violations is governed by a linear differential system involving the time-dependent NTK matrix $K(t)$. Following \cite{wang2022and}, this matrix admits the block representation
\begin{equation*}
	K(t)=
	\begin{bmatrix}
		K_{bb}(t) & K_{br}(t)\\
		K_{rb}(t) & K_{rr}(t)
	\end{bmatrix},
\end{equation*}
where $K_{bb}(t)$ and $K_{rr}(t)$ denote the NTK blocks associated with the boundary and PDE residual constraints, respectively, while $K_{br}(t)=K_{rb}^{T}(t)$ represents their cross-coupling. Under the standard assumptions of NTK theory, the limiting kernel $K(t)$ converges to a deterministic time-independent matrix $K$ in the infinite-width limit. Its spectral decomposition reveals that the components of the training error aligned with eigendirections of larger eigenvalues decay more rapidly, whereas those associated with smaller eigenvalues converge more slowly \cite{wang2022and,wang2021eigenvector}. Motivated by this spectral interpretation, the NTK-based scheme employs the traces of the diagonal kernel blocks as surrogate measures of their relative convergence rates, prescribing the candidate weights as
\begin{equation*}
	\begin{aligned}
		\hat{\lambda}_{bc}
		=
		\frac{\operatorname{tr}(K_{bb})+\operatorname{tr}(K_{rr})}
		{\operatorname{tr}(K_{bb})},\quad
		\hat{\lambda}_{r}
		=
		\frac{\operatorname{tr}(K_{bb})+\operatorname{tr}(K_{rr})}
		{\operatorname{tr}(K_{rr})}.
	\end{aligned}
\end{equation*}
For both weighting schemes, the loss weights are conventionally initialized to $\lambda_{bc}^{0}=\lambda_{r}^{0}=1$, with subsequent updates performed every $T$ training iterations. The comprehensive training guidelines established in \cite{wang2023expertsguidetrainingphysicsinformed} recommend setting the update interval to $T=1000$ and the exponential moving average factor to $\alpha=0.9$.

Beyond these loss-level weighting strategies, another line of research seeks to mitigate training imbalance through modifications to the optimizer itself. In the context of multi-task learning, the study in \cite{AdaTask_AAAI2023} exploits the linearity of Adam's first-moment recursion to distinguish the contributions of individual tasks to each parameter update. Based on these contributions, the degree of task dominance for each parameter is quantified through the relative proportions of exponentially averaged squared updates. The authors further link the observed task dominance to imbalanced gradient accumulation, whereby tasks with disproportionately large gradients can dominate the second-moment estimates used to determine adaptive learning rates, consequently exerting a disproportionate influence on parameter updates. To address this issue, AdaTask \cite{AdaTask_AAAI2023} applies Adam independently to each task's gradients and sums the resulting task-specific updates to obtain the overall parameter update. Similarly, MultiAdam \cite{pmlr-v202-yao23c} adopts a task-wise Adam optimization strategy for PINNs, but averages the resulting updates rather than summing them as in AdaTask.

Although these approaches can effectively promote balanced reductions across multiple loss components, existing explanations of why a particular task comes to dominate Adam's parameter updates rely either on GD-based theoretical formulations or on empirical characterizations of training behavior. To investigate this phenomenon theoretically, we draw inspiration from the relative accumulated update (rAU) indicator introduced in AdaTask and define $\Delta_k^{(i)}$ to quantify the contribution of the $i$-th task to the overall parameter update at iteration $k$. Using a Poisson problem as an illustrative example, we identify three characteristic phases in the evolution of $\{\Delta_k^{(i)}\}_{i=1}^{2}$. Focusing on the first phase, we derive an asymptotic relation for Adam configurations satisfying $(1-\beta_2)/(1-\beta_1)\ll1$, revealing how the shared second-moment estimate influences the step-to-step redistribution of task contributions through two distinct factors: disparities in curvature-related alignment and in the coordinatewise distributions of preconditioned task gradients.

\section{Methodology}
\subsection{Discrete-Time Analysis of Task Contribution Dynamics}
We begin by recalling the Adam iteration \cite{kingma2015iclr-adam} applied to the original unweighted PINN objective \eqref{pro:PINNs}:
\begin{align}
	m_k &= \beta_1 m_{k-1} + (1-\beta_1)g_{k-1}
	\label{eq:adam_m}\\
	v_k &= \beta_2 v_{k-1} + (1-\beta_2)g_{k-1}^2
	\label{eq:adam_v}\\
	\hat{m}_k &= \frac{m_k}{1-\beta_1^k},
	\hat{v}_k =\frac{v_k}{1-\beta_2^k}
	\label{eq:adam_bias}\\
	\theta_{k}
	&= \theta_{k-1}-\eta_k \frac{\hat{m}_k}{\sqrt{\hat{v}_k}+\delta}
	\label{eq:adam_theta}
\end{align}
Here, $\eta_k$ denotes the step size at iteration k, and $\delta>0$ is a small constant introduced to avoid division by zero. Throughout the subsequent analysis, vector powers, square roots, and quotients are understood componentwise, while $\odot$ denotes the Hadamard product. In the following analysis, we work in a fixed full-batch setting, so that $g_k$ is the full gradient of the problem \eqref{pro:PINNs}. For notational simplicity, we denote the gradients $\nabla_{\theta} \mathcal{L}_{bc}(\theta_k)$ and $\nabla_{\theta} \mathcal{L}_{r}(\theta_k)$ by $g^{(1)}_k$ and $g^{(2)}_k$, respectively. Thus,
\begin{equation*}
	g_k =g^{(1)}_k+g^{(2)}_k = \nabla_{\theta} \mathcal{L}_{bc}(\theta_k) + \nabla_{\theta} \mathcal{L}_{r}(\theta_k) = \nabla_{\theta} \mathcal{L}(\theta_k)
\end{equation*}
Exploiting the linear structure of the first-moment estimate \eqref{eq:adam_m}, \eqref{eq:adam_m}--\eqref{eq:adam_theta} can be equivalently expressed as \cite{AdaTask_AAAI2023,pmlr-v202-yao23c}
\begin{equation}
	\label{eq:split}
	\begin{aligned}
		m^{(i)}_k &= \beta_1 m^{(i)}_{k-1} + (1-\beta_1)g^{(i)}_{k-1} \; \; \  i=1,2\\
		v_k &= \beta_2 v_{k-1} + (1-\beta_2)(\sum\limits_{i=1}^{2}g^{(i)}_{k-1})^2\\
		\hat{v}_k & =\frac{v_k}{1-\beta_2^k},\quad \hat{m}^{(i)}_k  =\frac{m^{(i)}_k}{1-\beta_1^k} \; \; \  i=1,2\\
		\theta_{k}
		&= \theta_{k-1}-\eta_k \Big( \frac{\hat{m}^{(1)}_k}{\sqrt{\hat{v}_k}+\delta} + \frac{\hat{m}^{(2)}_k}{\sqrt{\hat{v}_k}+\delta} \Big)
	\end{aligned}
\end{equation}
Unlike vanilla gradient descent (GD), the influence of each loss term on the actual parameter update is no longer directly proportional to its raw gradient, but is encoded in its term-wise adaptive update component, defined as
\begin{equation*}
	\Delta \theta^{(i)}_k = \eta_k \frac{\hat{m}^{(i)}_k}{\sqrt{\hat{v}_k}+\delta}
\end{equation*}
Motivated by the relative accumulated update (rAU) metric introduced in AdaTask \cite{AdaTask_AAAI2023}, which characterizes task dominance at the level of individual parameters through exponentially averaged squared updates, we introduce a complementary measure that captures the relative magnitudes of task-wise update components across the entire shared parameter space at each iteration:
\begin{equation*}
	R^{(i)}_k=\frac{\Vert \Delta \theta^{(i)}_k \Vert}{\sum\limits_{j=1}^2\Vert \Delta \theta^{(j)}_k \Vert}
\end{equation*}
To facilitate the subsequent analysis, we assume that every component of the initial gradient is nonzero, ensuring that the second-moment estimates remain strictly positive componentwise throughout the iterations. Neglecting the small stabilization parameter $\delta$, the learning-rate and bias-correction factors then cancel in the normalized ratio, yielding the approximation
\begin{equation*}
	R^{(i)}_k \approx \frac{\Vert \frac{m^{(i)}_k}{\sqrt{v_k}} \Vert}{\sum\limits_{j=1}^{2} \Vert  \frac{m^{(j)}_k}{\sqrt{v_k}} \Vert} = \Delta^{(i)}_k
\end{equation*}

For illustration, we consider a one-dimensional Poisson problem.
\begin{equation*}
	\begin{aligned}
		&-\Delta u(x) = 4 \pi^2 \sin(2 \pi x) + 10000 \pi^2 \sin(100 \pi x),
		\quad x\in(0,1),\\
		&u(0)=u(1)=0
	\end{aligned}
\end{equation*}
The forcing term is chosen such that the exact solution is
\begin{equation*}
	u(x) = \sin(2 \pi x) + \sin(100 \pi x)
\end{equation*}
The solution is represented by a four-layer multilayer perceptron (MLP) incorporating the Fourier embedding proposed in \cite{wang2021eigenvector}, with the embedding scale set to $\sigma=13.5$. Training is conducted in a full-batch setting with
$N_{bc}=128$ and $N_r=2{,}000$. Optimization is performed using Adam with an initial learning rate of $10^{-3}$, $\beta_1=0.9$, and $\beta_2=0.999$, while the learning rate adopts an exponential decay schedule with a multiplicative factor of $0.9$ every $1000$ iterations.

\begin{figure}[htbp]
	\centering
	\begin{subfigure}[b]{0.40\textwidth}
		\centering
		\includegraphics[width=\textwidth]{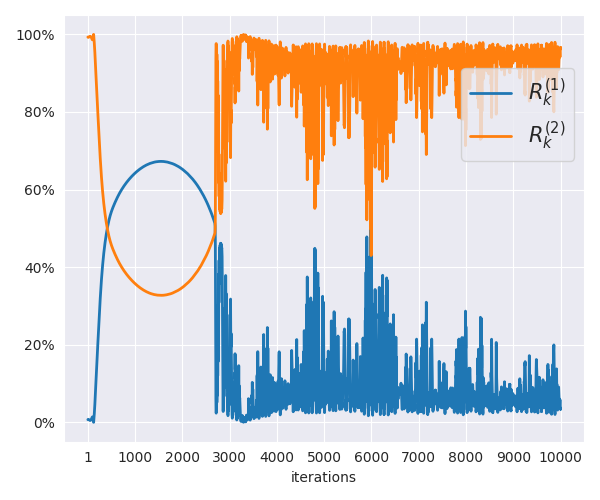}
		\caption{}
		\label{fig:R_ratio}
	\end{subfigure}
	\begin{subfigure}[b]{0.40\textwidth}
		\centering
		\includegraphics[width=\textwidth]{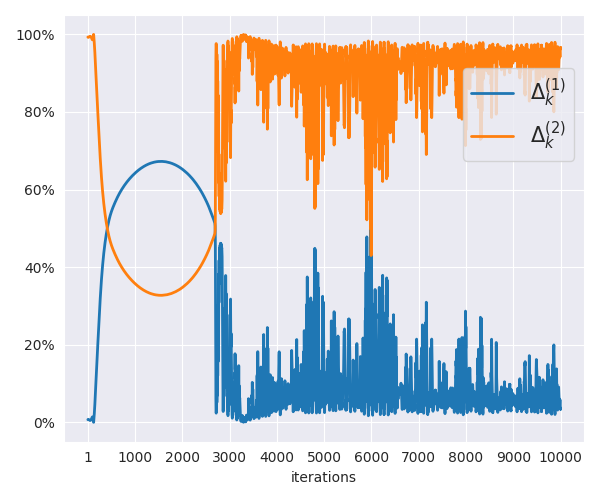}
		\caption{}
		\label{fig:Delta_ratio}
	\end{subfigure}
	\caption{Comparison of the actual and approximated update contributions:
		(a) $R_k^{(i)}$ and (b) $\Delta_k^{(i)}$.}
	\label{fig:update_contribution}
\end{figure}

As shown by the comparison between Fig.~\ref{fig:update_contribution}\subref{fig:R_ratio} and Fig.~\ref{fig:update_contribution}\subref{fig:Delta_ratio}, the approximated quantities closely track the actual contributions throughout the $10{,}000$ training iterations. This agreement supports the use of $\Delta_k^{(i)}$, $i=1,2$, as a reliable surrogate for the subsequent analysis.   
Based on the evolution exhibited in Fig.~\ref{fig:update_contribution}\subref{fig:Delta_ratio}, after excluding the first few transient iterations, we empirically identify three approximate dynamical regimes rather than imposing strict transition criteria: \textbf{Phase I (Smooth Evolution)}, \textbf{Phase II (Intermittent Spikes)} and \textbf{Phase III (Persistent Oscillations)}.
\begin{figure}[htbp]
	\centering
	\begin{subfigure}[b]{0.31\textwidth}
		\centering
		\includegraphics[width=\textwidth]{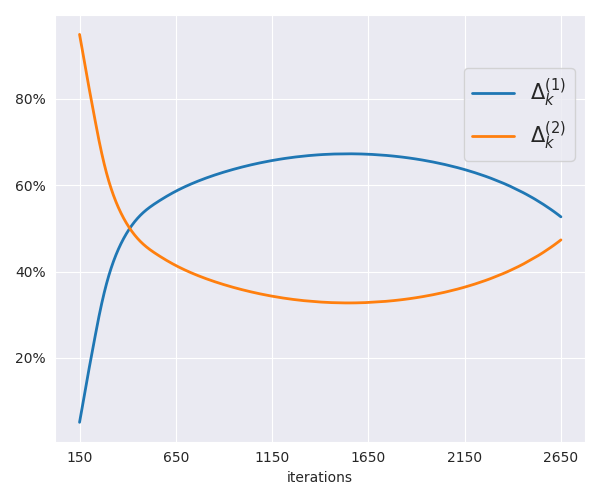}
		\caption{}
		\label{fig:phase1}
	\end{subfigure}
	\begin{subfigure}[b]{0.31\textwidth}
		\centering
		\includegraphics[width=\textwidth]{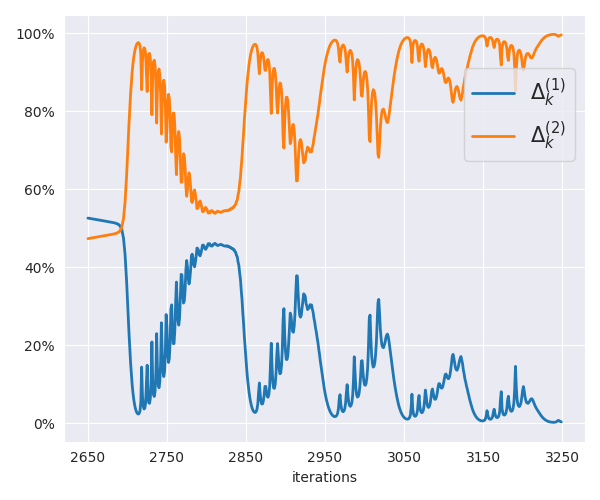}
		\caption{}
		\label{fig:phase2}
	\end{subfigure}
	\begin{subfigure}[b]{0.31\textwidth}
		\centering
		\includegraphics[width=\textwidth]{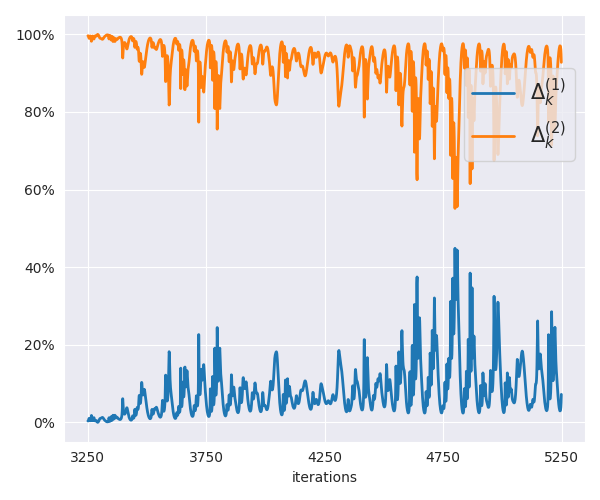}
		\caption{}
		\label{fig:phase3}
	\end{subfigure}
	\caption{Empirically identified dynamical regimes of the relative update
		contributions $\Delta_k^{(i)}$:
		(a) Phase I (Smooth Evolution);
		(b) Phase II (Intermittent Spikes);
		and (c) Phase III (Persistent Oscillations).}
	\label{fig:threephase}
\end{figure}

The salient dynamical features of these three regimes are displayed in Fig.~\ref{fig:threephase}. During Phase I [Fig.~\ref{fig:threephase}\subref{fig:phase1}], the two update contributions vary smoothly, with a gradual redistribution between the two loss components. The onset of Phase II [Fig.~\ref{fig:threephase}\subref{fig:phase2}] is marked by recurrent large-amplitude excursions, during which the relative contributions undergo rapid redistribution while repeatedly returning toward intermediate levels. By contrast, Phase III [Fig.~\ref{fig:threephase}\subref{fig:phase3}] is characterized by a persistent dominance pattern, where one contribution remains predominant and irregular oscillatory fluctuations persist around this strongly imbalanced state. The theoretical treatment developed below focuses on Phase I, with the aim of uncovering the mechanism by which Adam progressively reshapes the balance among loss-induced updates during the early training.

\begin{figure}[htbp]
	\centering
	\begin{subfigure}[b]{0.32\textwidth}
		\centering
		\includegraphics[width=\textwidth]{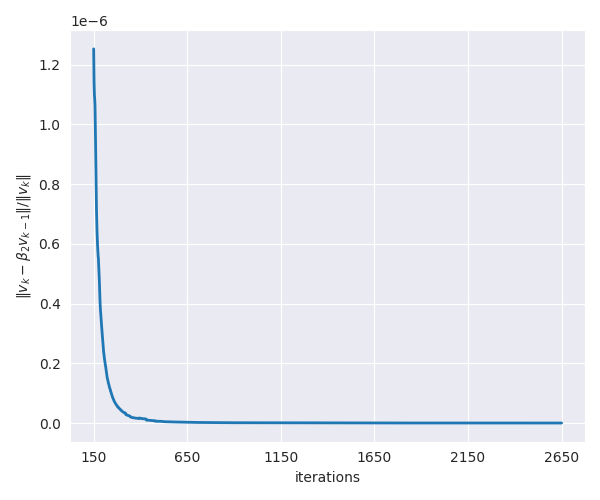}
		\caption{}
		\label{fig:second_moment_decay}
	\end{subfigure}
	\begin{subfigure}[b]{0.32\textwidth}
		\centering
		\includegraphics[width=\textwidth]{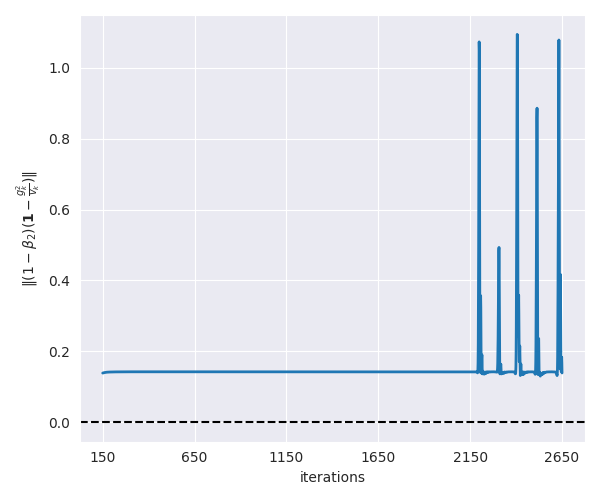}
		\caption{}
		\label{fig:q_norm}
	\end{subfigure}
	\begin{subfigure}[b]{0.32\textwidth}
		\centering
		\includegraphics[width=\textwidth]{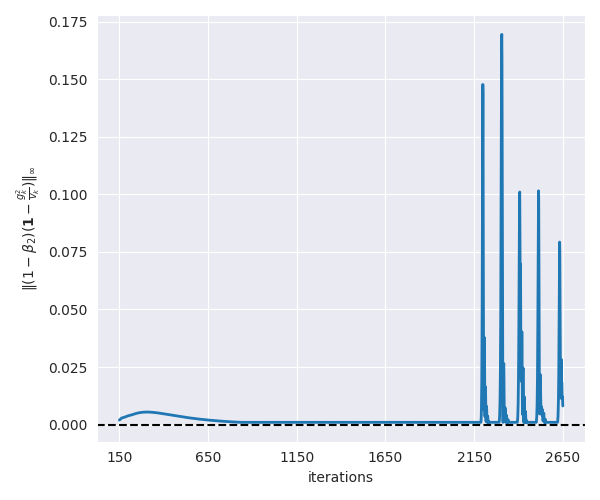}
		\caption{}
		\label{fig:q_infty}
	\end{subfigure}
	\caption{
		Numerical assessment of the second-moment approximation
		and first-order expansion.
		(a) Relative discrepancy
		$\Vert v_k-\beta_2v_{k-1}\Vert/\Vert v_k\Vert$.
		(b) Euclidean and (c) infinity norms of
		$(1-\beta_2)(\mathbf{1}-g_{k-1}^{2}/v_{k-1})$.
	}
	\label{fig:second_moment_expansion}
\end{figure}

The relative contributions of the two tasks can be characterized by the ratio $\rho_k:=\Delta_k^{(1)}/\Delta_k^{(2)}$. Motivated by the smooth step-to-step variation in this regime, we apply a first-order Taylor expansion of
$f(x,y)=\Vert x\Vert/\Vert y\Vert$ around the previous iterate, yielding $\rho_k \approx \rho_{k-1}\omega_k,$ where
\begin{equation*}
	\omega_k
	=
	1+
	\frac{
		\left\langle
		\frac{m_{k-1}^{(1)}}{\sqrt{v_{k-1}}},
		\frac{m_k^{(1)}}{\sqrt{v_k}}
		\right\rangle
	}{
		\left\|
		\frac{m_{k-1}^{(1)}}{\sqrt{v_{k-1}}}
		\right\|^2
	}
	-
	\frac{
		\left\langle
		\frac{m_{k-1}^{(2)}}{\sqrt{v_{k-1}}},
		\frac{m_k^{(2)}}{\sqrt{v_k}}
		\right\rangle
	}{
		\left\|
		\frac{m_{k-1}^{(2)}}{\sqrt{v_{k-1}}}
		\right\|^2
	}
\end{equation*}
Each inner-product ratio appearing in $\omega_k$ consists of a step-to-step magnitude ratio multiplied by a directional cosine. To further identify the factors governing the evolution of $\omega_k$, we substitute the task-wise first-moment recursion
into its definition. This decomposition separates the effects of historical momentum and newly incorporated task gradients, both modulated by the step-to-step variation in the shared second-moment estimate. Specifically, we obtain
\begin{equation*}
	\omega_k
	=
	1+\beta_1\mathcal{M}_k
	+(1-\beta_1)\mathcal{G}_k,
\end{equation*}
where
\begin{equation*}
	\begin{aligned}
		\mathcal{M}_k
		:={}&
		\frac{
			\left\langle
			\frac{m_{k-1}^{(1)}}{\sqrt{v_{k-1}}},
			\sqrt{\frac{v_{k-1}}{v_k}} \odot
			\frac{m_{k-1}^{(1)}}{\sqrt{v_{k-1}}}
			\right\rangle
		}{
			\left\|
			\frac{m_{k-1}^{(1)}}{\sqrt{v_{k-1}}}
			\right\|^2
		}
		-
		\frac{
			\left\langle
			\frac{m_{k-1}^{(2)}}{\sqrt{v_{k-1}}},
			\sqrt{\frac{v_{k-1}}{v_k}} \odot
			\frac{m_{k-1}^{(2)}}{\sqrt{v_{k-1}}}
			\right\rangle
		}{
			\left\|
			\frac{m_{k-1}^{(2)}}{\sqrt{v_{k-1}}}
			\right\|^2
		},
		\\[1.5ex]
		\mathcal{G}_k
		:={}&
		\frac{
			\left\langle
			\frac{m_{k-1}^{(1)}}{\sqrt{v_{k-1}}},
			\sqrt{\frac{v_{k-1}}{v_k}} \odot
			\frac{g_{k-1}^{(1)}}{\sqrt{v_{k-1}}}
			\right\rangle
		}{
			\left\|
			\frac{m_{k-1}^{(1)}}{\sqrt{v_{k-1}}}
			\right\|^2
		}
		-
		\frac{
			\left\langle
			\frac{m_{k-1}^{(2)}}{\sqrt{v_{k-1}}},
			\sqrt{\frac{v_{k-1}}{v_k}} \odot
			\frac{g_{k-1}^{(2)}}{\sqrt{v_{k-1}}}
			\right\rangle
		}{
			\left\|
			\frac{m_{k-1}^{(2)}}{\sqrt{v_{k-1}}}
			\right\|^2
		}.
	\end{aligned}
\end{equation*}
The componentwise scaling factor can be expressed explicitly in terms of the Adam second-moment recursion \eqref{eq:adam_v} as:
\[
\sqrt{\frac{v_{k-1}}{v_k}}
=
\left[
\mathbf{1}-(1-\beta_2)
\left(
\mathbf{1}-\frac{g_{k-1}^{2}}{v_{k-1}}
\right)
\right]^{(-1/2)}.
\]
Bai et al.~\cite{pmlr-v306-bai26o} characterized a decoupling regime in which the instantaneous squared-gradient contribution becomes negligible relative to the historical component, so that the second-moment estimate is predominantly governed by $\beta_2 v_{k-1}$. Consistent with this observation, Fig.\ref{fig:second_moment_expansion}\subref{fig:second_moment_decay} shows that the relative discrepancy $\Vert v_k-\beta_2v_{k-1} \Vert / \Vert v_k \Vert$ rapidly decreases to a near-zero level. Although the aggregate approximation may not hold uniformly across individual parameter coordinates, the near-geometric decay suggests a local expansion of the scaling factor. Specifically, for a sufficiently small expansion argument, the Maclaurin expansion yields
\[
\sqrt{\frac{v_{k-1}}{v_k}}
=
\mathbf{1}
+
\frac{1-\beta_2}{2}
\left(
\mathbf{1}-\frac{g_{k-1}^{2}}{v_{k-1}}
\right)
+\mathcal{O}\!\left(\| (1-\beta_2)\left(
\mathbf{1}-\frac{g_{k-1}^{2}}{v_{k-1}}
\right) \|^2\right).
\]
Figs.~\ref{fig:second_moment_expansion}\subref{fig:q_norm} and ~\ref{fig:second_moment_expansion}\subref{fig:q_infty} illustrate the evolution of the Euclidean and infinity norms of the expansion parameter, respectively. Both norms remain close to zero throughout much of the smooth-evolution regime, providing numerical support for the preceding first-order approximation.

The resulting first-order expression for $\omega_k$ involves only two sets of variables: $x_k^{(i)}=m_k^{(i)}/\sqrt{v_k}$ and $y_k^{(i)}=g_k^{(i)}/\sqrt{v_k}$. To further elucidate the factors governing its evolution, we establish an asymptotic relationship between these variables through a continuous-time analysis of Adam in the following subsection.

\subsection{Continuous-Time Analysis of Phase-I Dynamics}
While the discrete-time formulation accommodates iteration-dependent step sizes, we specialize to the constant-step setting $\eta_k=\eta>0$ for the continuous-time analysis, building on the Adam ODE framework established by Barakat and Bianchi \cite{barakat2021convergence}. Under suitable assumptions, they proved that the piecewise linear interpolations of stochastic Adam iterates converge in probability to the solution of the limiting ODE uniformly over compact time intervals, in the asymptotic regime characterized by
\[
\lim_{\eta\to0^+}\frac{1-\beta_j(\eta)}{\eta}=C_j>0,\qquad j=1,2.
\]
Building on this framework, Ma et al. \cite{pmlr-v145-ma22a} considered its full-batch specialization with the specific choice $\beta_j(\eta)=1-C_j\eta$, $j=1,2$, which is also adopted in this paper. For each $\eta>0$, let $z_i^{\eta}=(\theta_i^{\eta},m_i^{\eta},v_i^{\eta})$ denote the Adam state sequence generated with the hyperparameter tuple $(\eta,\beta_1(\eta),\beta_2(\eta))$, forming the family $\left\{\{z_i^{\eta}\}^{\infty}_{0}\right\}_{\eta>0}$. 
With $z^\eta(t)$ denoting the piecewise linear interpolation of $\{z_i^\eta\}_{i=0}^\infty$, the limiting dynamics take the form of the Adam ODE presented in Ma et al. \cite{pmlr-v145-ma22a}:
\begin{equation*}
	\begin{cases}
		\frac{d \theta(t)}{dt}
		=-\frac{(1-e^{-C_1t})^{-1}m(t)}
		{\sqrt{(1-e^{-C_2t})^{-1}v(t)}+\delta}\\
		\frac{d m(t)}{dt}
		=C_1\bigl(\nabla_{\theta}\mathcal{L}(\theta(t))-m(t)\bigr)\\
		\frac{d v(t)}{dt}
		=C_2\bigl(\nabla_{\theta}\mathcal{L}(\theta(t))^2-v(t)\bigr)\\
		\theta(0)=\theta_0,\; m(0)=v(0)=0.
	\end{cases}
\end{equation*}
Under the boundedness and global Lipschitz continuity of $\nabla\mathcal L$, the piecewise linear interpolations $z^\eta(t)$ converge uniformly on every compact time interval to $z(t)=(\theta(t),m(t),v(t))$ as $\eta\to0^+$. To bridge the preceding discrete-time analysis with the continuous-time formulation, we augment the Adam state sequence $z_i^\eta$ with the gradient variable $g^{\eta}_i$, yielding $\tilde{z}_i^\eta=(g^{\eta}_i,\theta_i^{\eta},m_i^{\eta},v_i^{\eta})$. The corresponding state dynamics are summarized in the following proposition.
\begin{proposition}
	\label{prop:adam_ode_limit}
	Let $\mathcal{L}:\mathbb{R}^d\to\mathbb{R}$ satisfy $\mathcal{L}\in C^2(\mathbb{R}^d)$, with $\nabla\mathcal{L}$ being Lipschitz continuous and bounded. For each $\eta>0$, consider the Adam state trajectory $\tilde{z}_i^\eta=(g_i^\eta,\theta_i^\eta,m_i^\eta,v_i^\eta)$ generated by \eqref{eq:adam_m}--\eqref{eq:adam_theta} under the parameterization $\beta_i(\eta)=1-C_i\eta$, $i=1,2$, with initial state $\tilde{z}_0^\eta=(g_0,\theta_0,0,0)$. Associated with this discrete sequence, define its piecewise linear interpolation $\tilde{z}^\eta:[0,\infty)\to\mathbb{R}^{4d}$ by
	\begin{equation*}
		\tilde{z}^\eta(t)
		=
		\tilde{z}_i^\eta
		+
		\frac{t-i\eta}{\eta}
		\left(
		\tilde{z}_{i+1}^\eta-\tilde{z}_i^\eta
		\right),
		\quad
		t\in[i\eta,(i+1)\eta), \, i\in\mathbb{N}_0.
	\end{equation*}
	Suppose that $\tilde{z}(t)=(g(t),\theta(t),m(t),v(t))$ denotes the solution of the following ODE system:
	\begin{equation}
		\label{eq:myODE}
		\begin{cases}
			\frac{d g(t)}{dt} = H(\theta(t)) \frac{d \theta(t)}{dt}\\
			\frac{d \theta(t)}{dt}
			=-\frac{(1-e^{-C_1t})^{-1}m(t)}
			{\sqrt{(1-e^{-C_2t})^{-1}v(t)}+\delta}\\
			\frac{d m(t)}{dt}
			=C_1\bigl(g(t)-m(t)\bigr)\\
			\frac{d v(t)}{dt}
			=C_2\bigl(g(t)^2-v(t)\bigr)\\
			g(0)=g_0, \; \theta(0)=\theta_0,\; m(0)=v(0)=0.
		\end{cases}
	\end{equation}
	Here, $H(\theta)=\nabla^2\mathcal{L}(\theta)$ is the Hessian matrix of $\mathcal{L}$ at $\theta$. Then, as $\eta\to0$, the family of interpolants ${\tilde{z}^{\eta}(t)}$ converges to the solution $\tilde{z}$ of \eqref{eq:myODE} uniformly on compact time intervals, i.e.,
	\begin{equation*}
		\lim\limits_{\eta \rightarrow 0}
		\sup_{t\in[0,T]}
		\Vert
		\tilde{z}^{\eta}(t)-\tilde{z}(t)
		\Vert=0,
		\quad \forall T>0.
	\end{equation*}
\end{proposition}
\begin{proof}
	Since $H(\theta)=\nabla^2\mathcal{L}(\theta)$, the chain rule yields
	\[
	\frac{d}{dt}\bigl(g(t)-\nabla\mathcal{L}(\theta(t))\bigr)=0,
	\]
	The initial condition $g_0=\nabla\mathcal{L}(\theta_0)$
	then implies $g(t)=\nabla\mathcal{L}(\theta(t))$
	for all $t\geq0$. Proposition 2 of \cite{pmlr-v145-ma22a} thus ensures that $\{ \theta^\eta \}_{\eta>0}$ converges to $\theta$ uniformly on compact time intervals as $\eta\to0^+$.
	For any $t \in [0,T]$, let $i$ be such that
	$t\in[i\eta,(i+1)\eta)$, and set
	$\alpha=(t-i\eta)/\eta$. The identity $g_i^\eta=\nabla\mathcal{L}(\theta_i^\eta)$, together with the $L$-Lipschitz continuity of $\nabla\mathcal{L}$, gives
	\[
	\begin{aligned}
		\|g^\eta(t)-g(t)\|
		&\leq
		\left\|(1-\alpha)\nabla\mathcal{L}(\theta_i^\eta)
		+\alpha\nabla\mathcal{L}(\theta_{i+1}^\eta)
		-\nabla\mathcal{L}(\theta^\eta(t))\right\|\\
		&\quad+\left\|\nabla\mathcal{L}(\theta^\eta(t))
		-\nabla\mathcal{L}(\theta(t))\right\|\\
		&\leq 2\alpha(1-\alpha)L
		\|\theta_{i+1}^\eta-\theta_i^\eta\|
		+L\|\theta^\eta(t)-\theta(t)\|
	\end{aligned}
	\]
	Moreover, the bias-corrected first moment satisfies
	\[
	\begin{aligned}
		\|\hat m_{i+1}^\eta\|
		\leq \max_{0\leq j\leq i}
		\|\nabla\mathcal{L}(\theta_j^\eta)\|
		\leq G
	\end{aligned}
	\]
	It follows that
	\[
	\begin{aligned}
		\|\theta_{i+1}^\eta-\theta_i^\eta\|
		=\eta\left\|
		\frac{\hat m_{i+1}^\eta}
		{\sqrt{\hat v_{i+1}^\eta}+\delta}
		\right\|
		\leq\frac{\eta G}{\delta}
	\end{aligned}
	\]
	Since $2\alpha(1-\alpha)\leq1/2$, combining the preceding estimates gives
	\[
	\begin{aligned}
		\sup_{t\in[0,T]}\|g^\eta(t)-g(t)\|
		\leq
		\frac{L\eta G}{2\delta}
		+
		L\sup_{t\in[0,T]}
		\|\theta^\eta(t)-\theta(t)\|
		\rightarrow0
		\quad (\eta\to0^+).
	\end{aligned}
	\]
	It remains to establish the convergence of the first- and
	second-moment estimates. For notational convenience, let
	$t_i=i\eta$ and $t_{i+1}=(i+1)\eta$. By the piecewise-linear interpolation and the integral
	representation of $m(t)$, we have
	\[
	\begin{aligned}
		\|m^\eta(t)-m(t)\|
		&=\Bigg\|
		m_i^\eta+\alpha(m_{i+1}^\eta-m_i^\eta)
		-m(t_i)
		-C_1\int_{t_i}^{t}(g(s)-m(s))\,ds
		\Bigg\|\\
		&\leq
		\|m_i^\eta-m(t_i)\|
		+\|m_{i+1}^\eta-m_i^\eta\| +
		C_1\int_{t_i}^{t_{i+1}}
		\|g(s)-m(s)\|\,ds.
	\end{aligned}
	\]
	The first-moment recursion and the uniform gradient bound
	yield
	\[
	\|m_{i+1}^\eta-m_i^\eta\|
	=C_1\eta\|g_i^\eta-m_i^\eta\|
	\leq 2C_1G\eta.
	\]
	Meanwhile, the variation-of-constants formula gives
	\[
	\begin{aligned}
		\|m(t)\|
		&=C_1\left\|
		\int_0^t e^{-C_1(t-s)}g(s)\,ds
		\right\| \leq G(1-e^{-C_1t})\leq G.
	\end{aligned}
	\]
	Consequently,
	\[
	\int_{t_i}^{t_{i+1}}
	\|g(s)-m(s)\|\,ds
	\leq 2G\eta.
	\]
	It therefore suffices to estimate the error at the grid points. Using the discrete first-moment recursion and the third equation of \eqref{eq:myODE}, we obtain
	\[
	\begin{aligned}
		\|m_i^\eta-m(t_i)\|
		=&\Bigg\|
		(1-C_1\eta)m_{i-1}^\eta
		+C_1\eta g_{i-1}^\eta
		-C_1\int_0^{t_i}(g(s)-m(s))\,ds
		\Bigg\|\\
		=&\Bigg\|
		(1-C_1\eta)
		\bigl(m_{i-1}^\eta-m(t_{i-1})\bigr)
		+C_1\eta
		\bigl(g_{i-1}^\eta-g(t_{i-1})\bigr) \\
		&+C_1\int_{t_{i-1}}^{t_i}
		\bigl[
		m(s)-m(t_{i-1})
		+g(t_{i-1})-g(s)
		\bigr]ds
		\Bigg\|
	\end{aligned}
	\]
	Since $\|m(t)\|\leq G$ , the continuous-time equations imply $\|\dot m(t)\|\leq 2C_1G$. Together with the $L$-Lipschitz continuity of
	$\nabla\mathcal L$, this gives, for $s\in[t_{i-1},t_i]$,
	\[
	\begin{aligned}
		\|m(s)-m(t_{i-1})\|
		&\leq 2C_1G(s-t_{i-1}),\\
		\|g(s)-g(t_{i-1})\|
		&\leq \frac{LG}{\delta}(s-t_{i-1}).
	\end{aligned}
	\]
	Hence, for $i \geq 1$,
	\[
	\begin{aligned}
		\|m_i^\eta-m(t_i)\|
		&\leq
		(1-C_1\eta)
		\|m_{i-1}^\eta-m(t_{i-1})\| + C_1\eta \sup_{t\in[0,T]}\|g^\eta(t)-g(t)\|
		+ C_1G\left(C_1+\frac{L}{2\delta}\right)\eta^2\\
		&\leq
		\Bigg(
		C_1\eta
		\sup_{t\in[0,T]}\|g^\eta(t)-g(t)\| + C_1G\left(C_1+\frac{L}{2\delta}\right)\eta^2
		\Bigg)
		\sum_{j=0}^{i-1}(1-C_1\eta)^j\\
		&\leq
		\sup_{t\in[0,T]}\|g^\eta(t)-g(t)\|
		+G\left(C_1+\frac{L}{2\delta}\right)\eta \rightarrow0
		\quad(\eta\to0^+).
	\end{aligned}
	\]
	Combining this estimate with the preceding interpolation
	bound, we conclude that
	\[
	\begin{aligned}
		\sup_{t\in[0,T]}\|m^\eta(t)-m(t)\|
		\leq
		\sup_{t\in[0,T]}\|g^\eta(t)-g(t)\|+
		G\left(5C_1+\frac{L}{2\delta}\right)\eta
		\rightarrow0
		\quad(\eta\to0^+).
	\end{aligned}
	\]
	The same argument applies to $v^\eta$, with details
	omitted. The convergence of all four components
	therefore establishes
	\[
	\lim_{\eta\to0^+}
	\sup_{t\in[0,T]}
	\|\tilde z^\eta(t)-\tilde z(t)\|=0,
	\]
	which completes the proof.
\end{proof}
The task-wise decomposition in \eqref{eq:split} carries over to the continuous-time limit, leading to the coupled ODE system below. The corresponding convergence follows by applying the preceding estimates componentwise, without requiring a separate proof.
\begin{equation}
	\label{eq:split_ODE}
	\begin{cases}
		\displaystyle
		\frac{d g^{(i)}(t)}{dt}
		=H^{(i)}(\theta(t))\frac{d\theta(t)}{dt},
		\quad i=1,2,\\[8pt]
		\displaystyle
		\frac{d\theta^{(i)}(t)}{dt}
		=-\frac{(1-e^{-C_1t})^{-1}
			\bigl(m^{(i)}(t)\bigr)}
		{\sqrt{(1-e^{-C_2t})^{-1}v(t)}+\delta},
		\quad i=1,2,
		\\[10pt]
		\displaystyle
		\frac{d m^{(i)}(t)}{dt}
		=C_1\bigl(g^{(i)}(t)-m^{(i)}(t)\bigr),
		\quad i=1,2,\\[8pt]
		\displaystyle
		\frac{d v(t)}{dt}
		=C_2\left[
		\bigl(g^{(1)}(t)+g^{(2)}(t)\bigr)^2-v(t)
		\right], \\[8pt]
		\displaystyle
		\theta(t) = \theta^{(1)}(t) + \theta^{(2)}(t), \\
		\displaystyle
		g^{(i)}(0)=g_0^{(i)},\quad
		m^{(i)}(0)=0,\quad i=1,2,\\
		\displaystyle
		\theta(0)=\theta_0,\quad v(0)=0.
	\end{cases}
\end{equation}
Here, $H^{(i)}(\theta)=\nabla_\theta^2\mathcal{L}_i(\theta)$ denotes the Hessian of the $i$-th loss component, $i=1,2$. We now consider the relationship between $x^{(i)}_{k}$ and $y^{(i)}_k$ from a continuous-time perspective. Under the initial-gradient condition imposed in Section~3.1, the second-moment equation in \eqref{eq:split_ODE}, together with $v(0)=0$, yields
\[
v_j(t)=C_2\int_0^t
e^{-C_2(t-s)}g_j(s)^2\,ds,
\quad j=1,\ldots,d.
\]
Since $g_{0,j}\neq 0$, continuity guarantees the existence of $a_j>0$ such that $|g_j(s)|\geq |g_{0,j}|/2$ for all $s\in[0,a_j]$. Consequently, for any $t>0$,
\[
v_j(t)
\geq
\frac{C_2g_{0,j}^2}{4}
\int_0^{\min\{t, a_j\}}
e^{-C_2(t-s)}\,ds
>0,
\quad j=1,\ldots,d.
\]
Then the componentwise positivity of $v(t)$ allows us to introduce the continuous-time counterparts of $x_k^{(i)}$ and $y_k^{(i)}$, together with an auxiliary variable $\bar z(t)$:
\[
x^{(i)}=\frac{m^{(i)}(t)}{\sqrt{v(t)}},
\qquad
y^{(i)}=\frac{g^{(i)}(t)}{\sqrt{v(t)}},
\qquad
\bar{z}=\frac{1}{C_2\sqrt{v(t)}},
\quad i=1,2.
\]
At the aggregate level, we set $x=x^{(1)}+x^{(2)}$ and $y=y^{(1)}+y^{(2)}$. 
Accordingly, we define the continuous-time counterparts of the indicators as
\begin{equation*}
	\begin{aligned}
		\mathcal{M}^{(i)}(t)
		=
		\frac{
			\left\langle
			x^{(i)}(t-\eta),
			\frac{\bar z(t)}{\bar z(t-\eta)}
			\odot x^{(i)}(t-\eta)
			\right\rangle
		}{
			\left\|x^{(i)}(t-\eta)\right\|^2
		}, \qquad
		\mathcal{G}^{(i)}(t)
		=
		\frac{
			\left\langle
			x^{(i)}(t-\eta),
			\frac{\bar z(t)}{\bar z(t-\eta)}
			\odot y^{(i)}(t-\eta)
			\right\rangle
		}{
			\left\|x^{(i)}(t-\eta)\right\|^2
		}
	\end{aligned}
\end{equation*}
for $i=1,2$. We further set
\[
\mathcal{M}(t)
=
\mathcal{M}^{(1)}(t)-\mathcal{M}^{(2)}(t),
\qquad
\mathcal{G}(t)
=
\mathcal{G}^{(1)}(t)-\mathcal{G}^{(2)}(t).
\]

Rescaling time as $\tau=C_2t$, we henceforth use an overdot to denote differentiation with respect to $\tau$. With the perturbation parameter $\varepsilon=C_2/C_1$, the resulting dynamics admit the following slow--fast formulation:
\begin{equation}
	\label{eq:split_slow_fast}
	\left\{
	\begin{aligned}
		\varepsilon\dot{x}^{(i)}
		&=y^{(i)}-x^{(i)}
		-\frac{\varepsilon}{2}
		(y^{2}-\mathbf{1})\odot x^{(i)},
		\quad i=1,2,\\
		\dot{y}^{(i)}
		&=-\frac12(y^{2}-\mathbf{1})
		\odot y^{(i)}
		-s_\varepsilon(\tau)
		\operatorname{diag}(\bar{z})
		H^{(i)}(\theta)P_\delta(\tau,\bar{z})x,
		\quad i=1,2,\\
		\dot{\bar{z}}
		&=-\frac12(y^{2}-\mathbf{1})\odot \bar{z},\\
		\dot{\theta}
		&=-\frac{s_\varepsilon(\tau)}{C_2}
		P_\delta(\tau,\bar{z})x,\\
		\dot{\tau} &=1.
	\end{aligned}
	\right.
\end{equation}
Here,
\[
s_\varepsilon(\tau)
=\frac{\sqrt{1-e^{-\tau}}}
{1-e^{-\tau/\varepsilon}},
\qquad
P_\delta(\tau,\bar{z})
=\operatorname{diag}\left(
\frac{1}{
	1+\delta C_2\sqrt{1-e^{-\tau}}\,\bar{z}}
\right).
\]
The small-parameter regime $\varepsilon\ll1$ is consistent with the standard Adam configuration $(\beta_1,\beta_2)=(0.9,0.999)$, which yields $\varepsilon=0.01$. Moreover, Zhang et al.~\cite{zhang2026adamconvergesmodificationupdate} establish convergence guarantees for sufficiently large $\beta_2$ under suitable conditions. For fixed $\beta_1<1$, increasing $\beta_2$ toward unity further drives $\varepsilon\to0^+$, motivating the following singular perturbation analysis. 
\begin{proposition}
	\label{prop:first_order_manifold}
	Let $\mathcal{L}_i\in C^4(\mathbb{R}^d)$, $i=1,2$, and consider the slow--fast system
	\eqref{eq:split_slow_fast} for $0<\varepsilon\ll1$. Setting $\varepsilon=0$, the critical manifold
	is given by
	\[
	\mathcal{C}_0=
	\left\{
	X=(x^{(1)},x^{(2)},y^{(1)},y^{(2)},
	\bar{z},\theta, \tau):
	x^{(i)}=y^{(i)},\ i=1,2
	\right\}.
	\]
	Let $\mathcal{S}_0$ be a compact submanifold of $\, \mathcal{C}_0$, possibly with boundary,
	contained in the region $\tau>\tau^*>0$. Then, for all sufficiently small $\varepsilon>0$, the following hold:
	
	\begin{enumerate}[
		label=(\roman*),
		labelindent=0pt,
		leftmargin=*
		]
		\item
		There exists a locally invariant,
		normally attracting slow manifold
		$\mathcal{S}_\varepsilon$
		diffeomorphic to $\mathcal{S}_0$,
		whose Hausdorff distance from
		$\mathcal{S}_0$ is $\mathcal{O}(\varepsilon)$.
		\item
		The slow manifold $\mathcal{S}_\varepsilon$
		admits the first-order approximation
		\begin{equation}
			\label{eq:first_order_manifold}
			x^{(i)}
			=y^{(i)}
			+\varepsilon s_0(\tau)
			\operatorname{diag}(\bar{z})
			H^{(i)}(\theta)
			P_\delta(\tau,\bar{z})y
			+\mathcal{O}(\varepsilon^2),
			\quad i=1,2,
		\end{equation}
		where $s_0(\tau)=\sqrt{1-e^{-\tau}}$.
	\end{enumerate}
\end{proposition}

\begin{proof}	
	We first establish (i). Introducing the fast time $\xi=\tau/\varepsilon$ and setting
	$\varepsilon=0$, the layer problem associated with \eqref{eq:split_slow_fast} takes the form
	\[
	\frac{dx^{(i)}}{d\xi}
	=y^{(i)}-x^{(i)},\qquad i=1,2,
	\]
	with the remaining variables held fixed. Hence, its equilibrium set is precisely $\mathcal{C}_0$.
	
	Since
	\[
	D_{(x^{(1)},x^{(2)})}
	\begin{pmatrix}
		y^{(1)}-x^{(1)}\\
		y^{(2)}-x^{(2)}
	\end{pmatrix}
	=-I_{2d},
	\]
	all eigenvalues are equal to $-1$. Thus, $\mathcal{C}_0$ is normally hyperbolic and attracting. By Fenichel's theorem \cite{FENICHEL197953,kuehn2015multiple}, for all sufficiently small $\varepsilon>0$, the compact submanifold $\mathcal{S}_0$ of $\mathcal{C}_0$ persists as a locally invariant, normally attracting slow manifold $\mathcal{S}_\varepsilon$ diffeomorphic to $\mathcal{S}_0$. Moreover, the Hausdorff distance between $\mathcal{S}_\varepsilon$ and $\mathcal{S}_0$ is $\mathcal{O}(\varepsilon)$, establishing (i).
	
	To prove (ii), we seek a local graph representation of $\mathcal{S}_\varepsilon$ in the form
	\begin{equation}
		\label{eq:slow_manifold_expansion_proof}
		x^{(i)}
		=h_\varepsilon^{(i)}
		=y^{(i)}
		+\varepsilon h_1^{(i)}
		+\mathcal{O}(\varepsilon^2),
		\qquad i=1,2.
	\end{equation}
	Here, $h_1^{(i)}$ depends on $(y^{(1)},y^{(2)},\bar z,\theta,\tau)$. Substituting \eqref{eq:slow_manifold_expansion_proof} into the first equation of
	\eqref{eq:split_slow_fast} yields
	\[
	\begin{aligned}
		\varepsilon\dot y^{(i)}
		+\mathcal{O}(\varepsilon^2)
		=
		-\varepsilon h_1^{(i)}
		-\frac{\varepsilon}{2}
		(y^{\odot2}-\mathbf{1})\odot y^{(i)}
		+\mathcal{O}(\varepsilon^2).
	\end{aligned}
	\]
	Since $\tau>\tau^*>0$, we have
	\[
	\begin{aligned}
		s_\varepsilon(\tau)-s_0(\tau)
		=
		s_0(\tau)
		\frac{e^{-\tau/\varepsilon}}
		{1-e^{-\tau/\varepsilon}} =\mathcal{O} \left(e^{-\tau^*/\varepsilon}\right) =o(\varepsilon).
	\end{aligned}
	\]
	Thus, $s_\varepsilon(\tau)=s_0(\tau)+o(\varepsilon)$ uniformly for $\tau>\tau^*$. Similarly, substituting the expansion \eqref{eq:slow_manifold_expansion_proof} into the second equation gives
	\[
	\begin{aligned}
		\dot y^{(i)}
		=
		-\frac12(y^{2}-\mathbf{1})
		\odot y^{(i)}
		-s_0(\tau)
		\operatorname{diag}(\bar z)
		H^{(i)}(\theta)
		P_\delta(\tau,\bar z)y
		+\mathcal{O}(\varepsilon).
	\end{aligned}
	\]
	Combining these expressions and comparing the coefficients of $\varepsilon$, we obtain
	\[
	h_1^{(i)}
	=s_0(\tau)
	\operatorname{diag}(\bar z)
	H^{(i)}(\theta)
	P_\delta(\tau,\bar z)y,
	\qquad i=1,2.
	\]
	Substituting this into \eqref{eq:slow_manifold_expansion_proof} establishes (ii).
\end{proof}

By the local invariance of $\mathcal{S}_\varepsilon$, trajectories starting on this manifold remain on it until they reach its boundary. Moreover, its normal attraction guarantees the existence of a neighborhood of $\mathcal{S}_\varepsilon$ such that trajectories whose initial states lie in this neighborhood approach the manifold exponentially fast. These properties provide the theoretical basis for applying the first-order slow-manifold approximation to the continuous-time factors $\mathcal{M}(t)$ and $\mathcal{G}(t)$.
\begin{proposition}
	\label{prop:mu_kappa_manifold11}
	Under the assumptions of Proposition~\ref{prop:first_order_manifold}, consider a trajectory whose initial state is sufficiently close to $\mathcal{S}_\varepsilon$,
	with $\|y^{(i)}\|$ uniformly bounded away from zero, $i=1,2$. Let $t_k=k\eta$ and define
	\[
	A^{(i)}(t):=
	s_0(C_2t) \operatorname{diag}(\bar z(t))
	H^{(i)}(\theta(t))
	P_\delta(C_2t,\bar z(t)),
	\qquad i=1,2.
	\]
	Then, for $k>1$, the following hold:
	\begin{enumerate}[
		label=(\roman*),
		labelindent=0pt,
		leftmargin=*
		]
		\item
		For $i=1,2$, define
		\[
		\mu^{(i)}(t)
		:=
		\frac{\|y^{(i)}(t)\|}
		{\|x^{(i)}(t)\|}.
		\]
		Then $\mu^{(i)}(t)$ admits the first-order
		approximation
		\begin{equation}
			\label{eq:mu_manifold_approximation}
			\mu^{(i)}(t)
			=
			1-\varepsilon \frac{
				\left\langle
				y^{(i)}(t),
				A^{(i)}(t)y(t)
				\right\rangle
			}{
				\|y^{(i)}(t)\|^2
			}
			+\mathcal{O}(\varepsilon^2).
		\end{equation}
		\item
		$\mathcal{M}(t)$ and $\mathcal{G}(t)$ admit the following approximations:
		\begin{equation}
			\label{eq:M_G_manifold_relation}
			\begin{aligned}
				\mathcal{M}(t_{k+1})
				={}&
				\frac{C_2\eta}{2}
				\Bigg[
				\frac{
					\|y(t_k)\odot y^{(2)}(t_k)\|^2
				}{
					\|y^{(2)}(t_k)\|^2
				}
				-
				\frac{
					\|y(t_k)\odot y^{(1)}(t_k)\|^2
				}{
					\|y^{(1)}(t_k)\|^2
				}
				\Bigg]+\mathcal{O}\left(
				C_2\eta\varepsilon+C_2^2\eta^2
				\right),
				\\
				\mathcal{G}(t_{k+1})
				={}&
				\varepsilon 
				\Bigg[
				\frac{
					\left\langle
					y^{(2)}(t_k),
					A^{(2)}(t_k)y(t_k)
					\right\rangle
				}{
					\|y^{(2)}(t_k)\|^2
				}
				-
				\frac{
					\left\langle
					y^{(1)}(t_k),
					A^{(1)}(t_k)y(t_k)
					\right\rangle
				}{
					\|y^{(1)}(t_k)\|^2
				}
				\Bigg]\\
				&+
				\frac{C_2\eta}{2}
				\Bigg[
				\frac{
					\|y(t_k)\odot y^{(2)}(t_k)\|^2
				}{
					\|y^{(2)}(t_k)\|^2
				}
				-
				\frac{
					\|y(t_k)\odot y^{(1)}(t_k)\|^2
				}{
					\|y^{(1)}(t_k)\|^2
				}
				\Bigg]
				+\mathcal{O}\left(
				\varepsilon^2
				+C_2\eta\varepsilon
				+C_2^2\eta^2
				\right).
			\end{aligned}
		\end{equation}
	\end{enumerate}
\end{proposition}
\begin{proof}
	Let
	\[
	B(t)=\frac12\operatorname{diag}
	\left(\mathbf{1}-y^{2}(t)\right).
	\]
	Returning to the original time variable $t$ via $\tau=C_2t$, Proposition~\ref{prop:first_order_manifold} yields
	\begin{equation}
		\label{eq:proof_manifold}
		x^{(i)}(t)
		=
		y^{(i)}(t)
		+\varepsilon A^{(i)}(t)y(t)
		+\mathcal{O}(\varepsilon^2),
		\qquad i=1,2.
	\end{equation}
	We first prove (i). Applying the first-order Taylor expansion of the Euclidean norm to
	\eqref{eq:proof_manifold}, we obtain
	\begin{equation}
		\label{eq:proof_mu_expansion}
		\begin{aligned}
			\mu^{(i)}(t)
			&=
			\frac{
				\|y^{(i)}(t)\|
			}{
				\|x^{(i)}(t)\|
			}=
			\frac{
				\|y^{(i)}(t)\|
			}{
				\|y^{(i)}(t)+
				\varepsilon A^{(i)}(t)y(t)
				+\mathcal{O}(\varepsilon^2)\|
			}\\
			&=
			1-\varepsilon
			\frac{
				\left\langle
				y^{(i)}(t),
				A^{(i)}(t)y(t)
				\right\rangle
			}{
				\|y^{(i)}(t)\|^2
			}
			+\mathcal{O}(\varepsilon^2).
		\end{aligned}
	\end{equation}
	This establishes (i). We now prove (ii). In the original time scale, the evolution equation for $\bar z$ in \eqref{eq:split_slow_fast} becomes
	\[
	\frac{d\bar z}{dt}=C_2B(t)\bar z(t).
	\]
	Applying a first-order Taylor expansion at $t_k$ gives
	\begin{equation}
		\label{eq:proof_z_ratio}
		\frac{\bar z(t_{k+1})}{\bar z(t_k)}
		=
		\mathbf{1}
		+C_2\eta B(t_k)\mathbf{1}
		+\mathcal{O}(C_2^2\eta^2).
	\end{equation}
	We first consider $\mathcal{M}^{(i)}(t_{k+1})$. Substituting \eqref{eq:proof_z_ratio}
	into its definition yields
	\begin{equation}
		\label{eq:proof_M_expansion}
		\begin{aligned}
			\mathcal{M}^{(i)}(t_{k+1})
			&=
			\frac{
				\left\langle
				x^{(i)}(t_k),
				\dfrac{\bar z(t_{k+1})}{\bar z(t_k)}
				\odot x^{(i)}(t_k)
				\right\rangle
			}{
				\|x^{(i)}(t_k)\|^2
			}=
			1+C_2\eta
			\frac{
				\left\langle
				x^{(i)}(t_k),
				B(t_k)x^{(i)}(t_k)
				\right\rangle
			}{
				\|x^{(i)}(t_k)\|^2
			}
			+\mathcal{O}(C_2^2\eta^2)\\
			&=
			1+C_2\eta
			\frac{
				\left\langle
				y^{(i)}(t_k),
				B(t_k)y^{(i)}(t_k)
				\right\rangle
			}{
				\|y^{(i)}(t_k)\|^2
			}+
			\mathcal{O}\left(
			C_2\eta\varepsilon+C_2^2\eta^2
			\right),
		\end{aligned}
	\end{equation}
	where the last equality follows from
	\eqref{eq:proof_manifold}.
	By the definition of $B(t)$,
	\begin{equation}
		\label{eq:proof_B_identity}
		\begin{aligned}
			\frac{
				\left\langle
				y^{(i)}(t_k),
				B(t_k)y^{(i)}(t_k)
				\right\rangle
			}{
				\|y^{(i)}(t_k)\|^2
			}=
			\frac12\left[
			1-
			\frac{
				\|y(t_k)\odot y^{(i)}(t_k)\|^2
			}{
				\|y^{(i)}(t_k)\|^2
			}
			\right].
		\end{aligned}
	\end{equation}
	Taking the difference between \eqref{eq:proof_M_expansion} for $i=1$ and $i=2$ and applying \eqref{eq:proof_B_identity} establishes the approximation for $\mathcal{M}(t_{k+1})$.
	
	We next consider $\mathcal{G}^{(i)}(t_{k+1})$. Combining \eqref{eq:proof_manifold}
	and \eqref{eq:proof_z_ratio}, we obtain
	\begin{equation}
		\label{eq:proof_G_expansion}
		\begin{aligned}
			\mathcal{G}^{(i)}(t_{k+1})
			&=
			\frac{
				\left\langle
				x^{(i)}(t_k),
				\dfrac{\bar z(t_{k+1})}{\bar z(t_k)}
				\odot y^{(i)}(t_k)
				\right\rangle
			}{
				\|x^{(i)}(t_k)\|^2
			}\\
			&=
			\frac{
				\left\langle
				y^{(i)}(t_k)
				+\varepsilon A^{(i)}(t_k)y(t_k)
				+\mathcal{O}(\varepsilon^2),
				y^{(i)}(t_k)
				+C_2\eta B(t_k)y^{(i)}(t_k)
				+\mathcal{O}(C_2^2\eta^2)
				\right\rangle
			}{
				\left\|
				y^{(i)}(t_k)
				+\varepsilon A^{(i)}(t_k)y(t_k)
				+\mathcal{O}(\varepsilon^2)
				\right\|^2
			}\\
			&=
			1-\varepsilon
			\frac{
				\left\langle
				y^{(i)}(t_k),
				A^{(i)}(t_k)y(t_k)
				\right\rangle
			}{
				\|y^{(i)}(t_k)\|^2
			}+
			C_2\eta
			\frac{
				\left\langle
				y^{(i)}(t_k),
				B(t_k)y^{(i)}(t_k)
				\right\rangle
			}{
				\|y^{(i)}(t_k)\|^2
			}+
			\mathcal{O}\left(
			\varepsilon^2
			+C_2\eta\varepsilon
			+C_2^2\eta^2
			\right).
		\end{aligned}
	\end{equation}
	Finally, taking the difference between \eqref{eq:proof_G_expansion} for $i=1$ and $i=2$ and using \eqref{eq:proof_B_identity} gives the claimed approximation for
	$\mathcal{G}(t_{k+1})$, completing the proof.
\end{proof}

\subsection{Shared Preconditioning and Task-Contribution Imbalance in Adam}
\begin{figure}[htbp]
	\centering
	\begin{subfigure}[b]{0.32\textwidth}
		\centering
		\includegraphics[width=\textwidth]{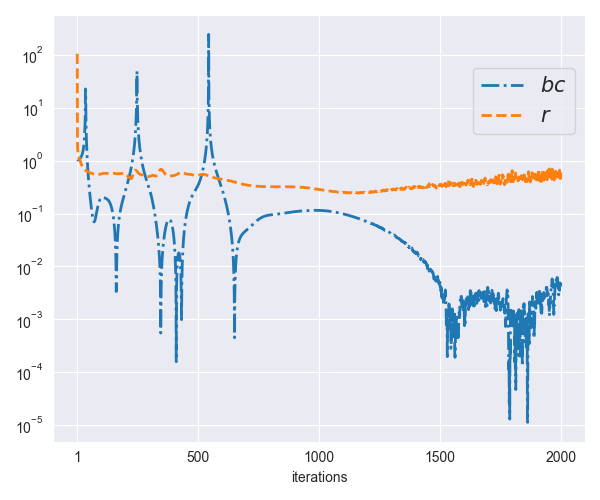}
		\caption{}
		\label{fig:manifold_lr_1}
	\end{subfigure}
	\begin{subfigure}[b]{0.32\textwidth}
		\centering
		\includegraphics[width=\textwidth]{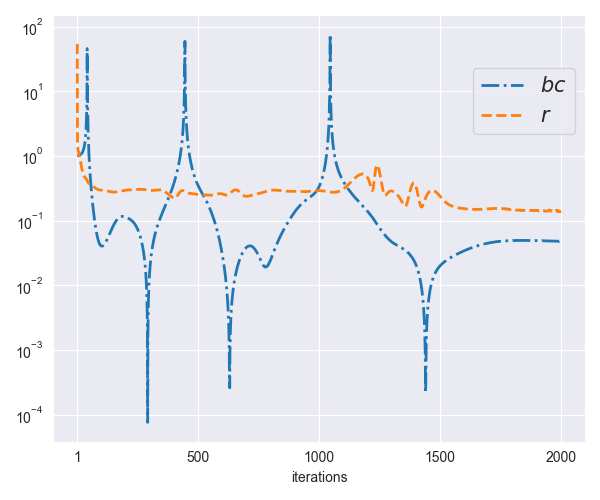}
		\caption{}
		\label{fig:manifold_lr_2}
	\end{subfigure}
	\begin{subfigure}[b]{0.32\textwidth}
		\centering
		\includegraphics[width=\textwidth]{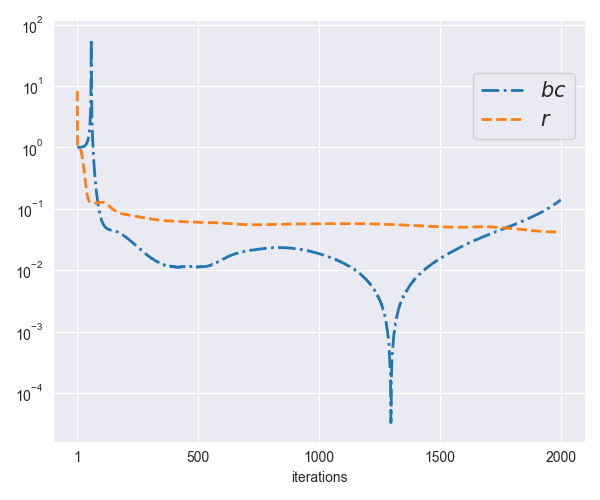}
		\caption{}
		\label{fig:manifold_lr_3}
	\end{subfigure}
	\caption{
		Relative $L^2$ errors of the first-order
		predictions of $x_k^{(i)}-y_k^{(i)}$
		given by \eqref{eq:discrete_slow_manifold},
		with respect to their actual values from
		the discrete Adam iterations in the 1D Poisson problem.
		For $(\beta_1,\beta_2)=(0.9,0.999)$:
		(a) $\eta=10^{-4}$;
		(b) $\eta=5\times10^{-5}$;
		(c) $\eta=10^{-5}$.
	}
	\label{fig:manifold_learning_rate}
\end{figure}

\begin{figure}[htbp]
	\centering
	\begin{subfigure}[b]{0.32\textwidth}
		\centering
		\includegraphics[width=\textwidth]{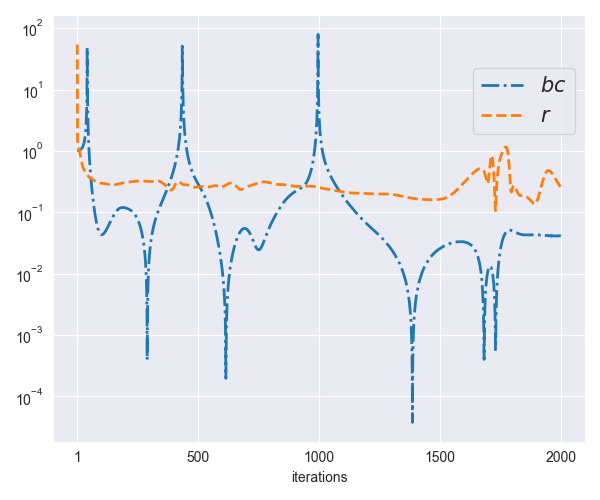}
		\caption{}
		\label{fig:manifold_beta_1}
	\end{subfigure}
	\begin{subfigure}[b]{0.32\textwidth}
		\centering
		\includegraphics[width=\textwidth]{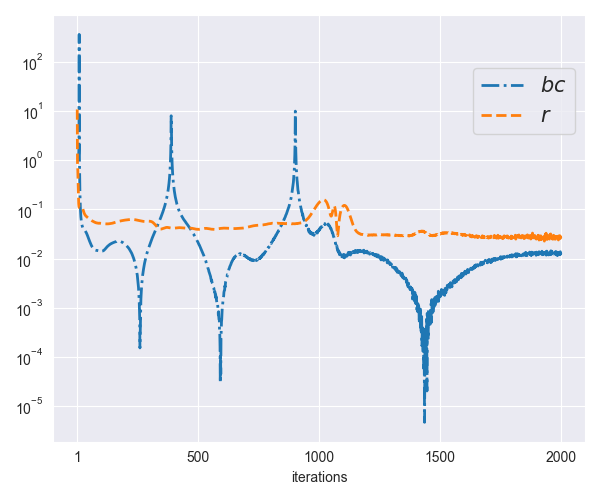}
		\caption{}
		\label{fig:manifold_beta_2}
	\end{subfigure}
	\begin{subfigure}[b]{0.32\textwidth}
		\centering
		\includegraphics[width=\textwidth]{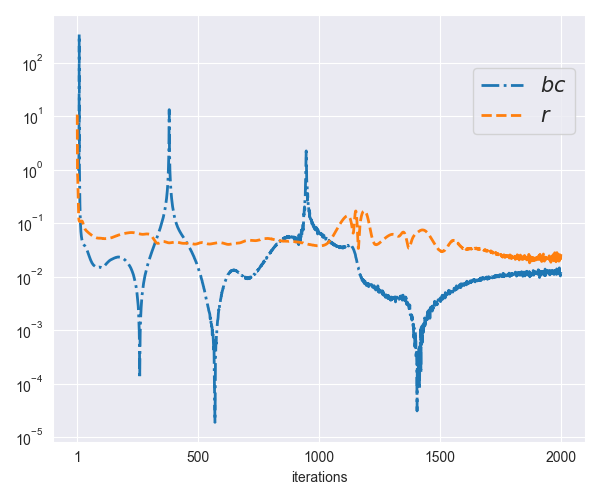}
		\caption{}
		\label{fig:manifold_beta_3}
	\end{subfigure}
	\caption{
		Relative $L^2$ errors of the first-order
		predictions of $x_k^{(i)}-y_k^{(i)}$
		given by \eqref{eq:discrete_slow_manifold},
		with respect to their actual values from
		the discrete Adam iterations in the 1D Poisson problem.
		For $\eta=5\times10^{-5}$:
		(a) $(\beta_1,\beta_2)=(0.9,0.99995)$;
		(b) $(\beta_1,\beta_2)=(0.5,0.999)$;
		(c) $(\beta_1,\beta_2)=(0.5,0.99995)$.
	}
	\label{fig:manifold_adam_parameters}
\end{figure}

Proposition~\ref{prop:adam_ode_limit} establishes a continuous-time asymptotic characterization of the discrete Adam dynamics. For sufficiently small $\eta$, provided that the associated continuous-time trajectory enters a sufficiently small neighborhood of the attracting slow manifold $\mathcal{S}_\varepsilon$, the first-order relation \eqref{eq:first_order_manifold}, evaluated at $t_k=k\eta$, suggests the following discrete approximation:
\begin{equation}
	\label{eq:discrete_slow_manifold}
	\begin{aligned}
		x_k^{(i)}
		\approx{}&
		y_k^{(i)}
		+\varepsilon s_{0,k}
		\frac{1}{C_2\sqrt{v_k}}
		\odot H^{(i)}(\theta_k)
		\left(
		\frac{g_k}
		{\sqrt{v_k}+\delta s_{0,k}}
		\right),
		\quad i=1,2,
	\end{aligned}
\end{equation}
where $s_{0,k}=\sqrt{1-e^{-C_2k\eta}}$. Figures~\ref{fig:manifold_learning_rate}
and~\ref{fig:manifold_adam_parameters} assess the accuracy of the discrete approximation
\eqref{eq:discrete_slow_manifold} under different learning rates and Adam hyperparameters, respectively. For fixed $(\beta_1,\beta_2)=(0.9,0.999)$,
Fig.~\ref{fig:manifold_learning_rate} shows that the relative errors exhibit transient spikes but remain relatively small over substantial portions of the subsequent iterations. Reducing the learning rate generally lowers the errors for both loss components, while also suppressing pronounced error spikes. At a fixed learning rate $\eta=5\times10^{-5}$, Fig.~\ref{fig:manifold_adam_parameters} examines the effect of the singular perturbation parameter $\varepsilon=(1-\beta_2)/(1-\beta_1)$. Comparing panels (a) and (c), reducing $\varepsilon$ from $5\times10^{-4}$ to $10^{-4}$ improves the
approximation accuracy for both components, while panels (b) and (c) exhibit comparably small errors despite their different values of $\varepsilon$. These findings provide further numerical support for the applicability of the first-order slow-manifold approximation to the discrete Adam dynamics under the considered configurations.

Building on this evidence, we use the relations established in Proposition~\ref{prop:mu_kappa_manifold11} to reformulate $\omega_k$, thereby explicitly separating the terms governing the relative update contributions of the two tasks. To this end, we first introduce the following task-specific quantities for $i=1,2$:
\begin{equation*}
	\begin{aligned}
		\kappa^{(i)}(t)
		:=
		\frac{
			\left\langle
			y^{(i)}(t),
			\displaystyle
			A^{(i)}(t) y(t)
			\right\rangle
		}{
			\|y^{(i)}(t)\|^2
		},\qquad
		\chi^{(i)}(t)
		:=
		\frac{
			\|y(t)\odot y^{(i)}(t)\|^2
		}{
			\|y^{(i)}(t)\|^2
		}.
	\end{aligned}
\end{equation*}
With these definitions, the continuous-time analogue of $\omega_k$, evaluated at $t_k=k\eta$, can be expressed as
\begin{equation*}
	\label{eq:omega_continuous}
	\begin{aligned}
		\omega(t_k)
		\approx{}
		1+(1-\beta_1)\varepsilon
		\underbrace{
			\left[
			\kappa^{(2)}(t_{k-1})
			-\kappa^{(1)}(t_{k-1})
			\right]
		}_{\text{term A}}+\frac{1-\beta_2}{2}
		\underbrace{
			\left[
			\chi^{(2)}(t_{k-1})
			-\chi^{(1)}(t_{k-1})
			\right]
		}_{\text{term B}}.
	\end{aligned}
\end{equation*}
The term A captures the difference in curvature-related alignment between the two tasks. Specifically, $\kappa^{(i)}(t)$ reflects the projection of $A^{(i)}(t)y(t)$ onto the direction of $y^{(i)}(t)$, scaled by the magnitude of the latter. The term B, by contrast, admits the following coordinatewise decomposition:
\begin{equation*}
	\begin{aligned}
		\chi^{(2)}(t)-\chi^{(1)}(t)=
		\sum_{j=1}^{d} y_j^2(t)
		\left(
		\frac{(y_j^{(2)}(t))^2}
		{\|y^{(2)}(t)\|^2}
		-
		\frac{(y_j^{(1)}(t))^2}
		{\|y^{(1)}(t)\|^2}
		\right).
	\end{aligned}
\end{equation*}
The difference in parentheses reflects the disparity in the concentration of the magnitudes of $y_{k-1}^{(1)}$ and $y_{k-1}^{(2)}$ at a given parameter coordinate. This disparity is amplified in the weighted sum when the corresponding entries of the two loss gradients have the same sign. If the two loss gradients exhibit identical relative distributions of magnitudes across all parameter coordinates, the term B necessarily vanishes. Such an exact agreement, however, cannot be assumed in general, as Wang et al.~\cite{wang2021understanding} documented pronounced disparities in the distributions of back-propagated gradients associated with boundary and PDE-residual losses. 

Both terms depend on Adam's shared second-moment estimate $v(t)$, highlighting two distinct factors through which the shared preconditioner affects the evolution of relative task update contributions. When this estimate is dominated by one task, the resulting preconditioner may exert asymmetric influence on the two tasks, potentially limiting Adam's ability to redistribute their update contributions.

\begin{figure}[H]
	\centering
	\begin{subfigure}[b]{0.32\textwidth}
		\centering
		\includegraphics[width=\linewidth]{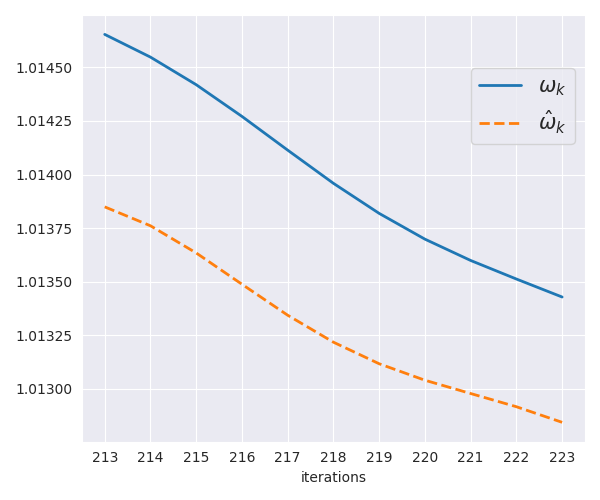}
		\caption{}
		\label{fig:omega_hat_early}
	\end{subfigure}
	\begin{subfigure}[b]{0.32\textwidth}
		\centering
		\includegraphics[width=\linewidth]{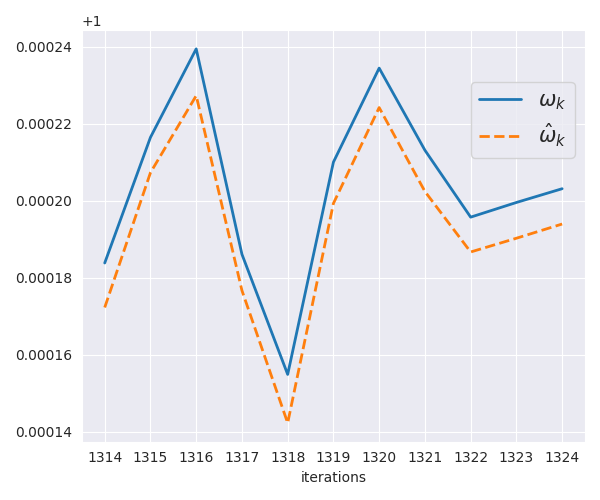}
		\caption{}
		\label{fig:omega_hat_intermediate}
	\end{subfigure}
	\begin{subfigure}[b]{0.32\textwidth}
		\centering
		\includegraphics[width=\linewidth]{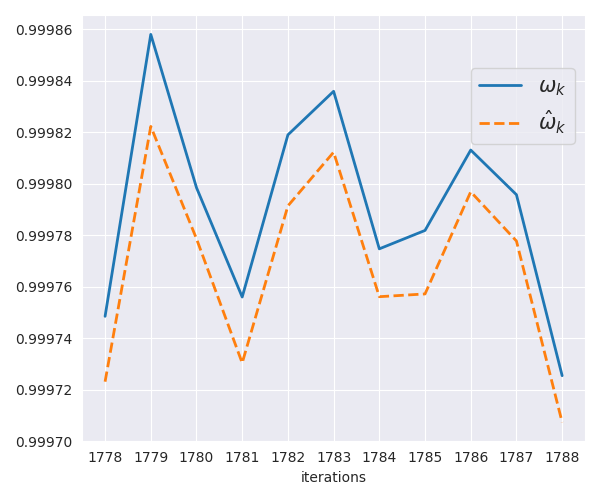}
		\caption{}
		\label{fig:omega_hat_later}
	\end{subfigure}
	\caption{
		Comparison of $\omega_k$ and $\hat{\omega}_k$
		for the one-dimensional Poisson problem
		under an exponentially decaying learning rate,
		over three randomly selected iteration intervals
		within \textbf{Phase I}:
		(a) $k=213$--$223$;
		(b) $k=1314$--$1324$;
		(c) $k=1778$--$1788$.
		Different vertical scales are used to
		resolve the local variations.
	}
	\label{fig:omega_hat_comparison}
\end{figure}

\begin{figure}[htbp]
	\centering
	\begin{subfigure}[b]{0.32\textwidth}
		\centering
		\includegraphics[width=\linewidth]{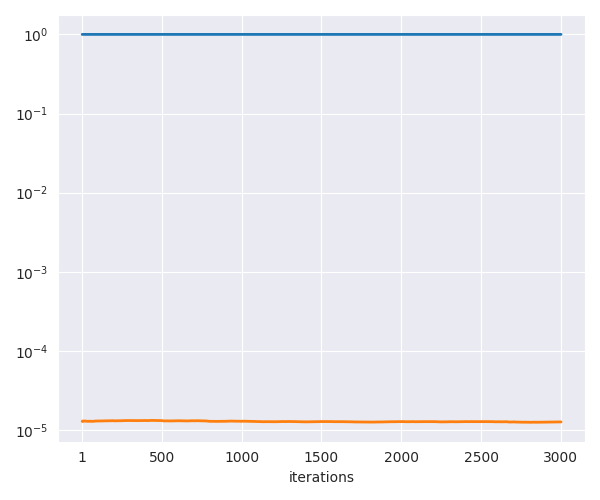}
		\caption{}
		\label{fig:secondmoment_adam}
	\end{subfigure}
	\begin{subfigure}[b]{0.32\textwidth}
		\centering
		\includegraphics[width=\linewidth]{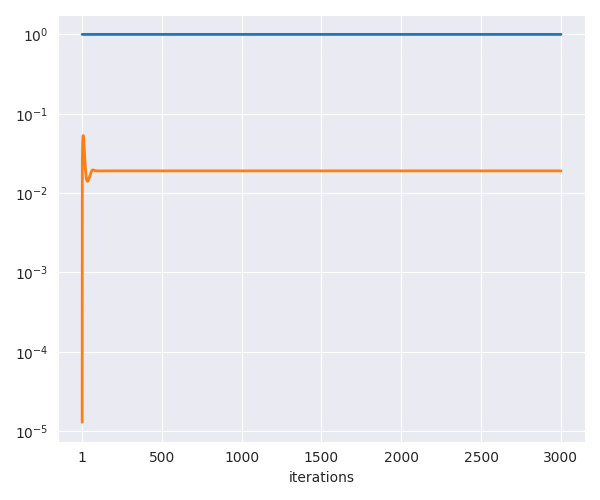}
		\caption{}
		\label{fig:secondmoment_gradnorm}
	\end{subfigure}
	\begin{subfigure}[b]{0.32\textwidth}
		\centering
		\includegraphics[width=\linewidth]{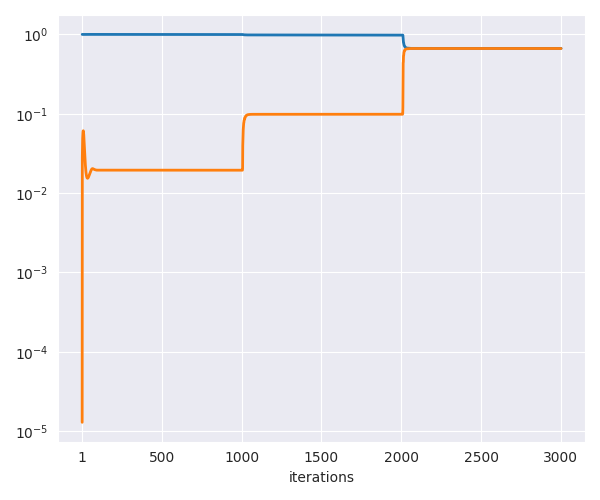}
		\caption{}
		\label{fig:secondmoment_ntk}
	\end{subfigure}
	\caption{
		Relative $L^2$ errors of the task-specific
		second-moment estimates $v^{(1)}$
		and $v^{(2)}$ with respect to the shared
		estimate $v$, for Adam trained
		(a) without loss reweighting,
		(b) with gradient-norm-based loss reweighting,
		and (c) with NTK-based loss reweighting.
		The blue and orange curves represent
		$\|v-v^{(1)}\|_2/\|v\|_2$
		and $\|v-v^{(2)}\|_2/\|v\|_2$,
		respectively.
		The vertical axes use a logarithmic scale.
	}
	
	\label{fig:secondmoment_comparison}
\end{figure}

Moreover, Proposition~\ref{prop:mu_kappa_manifold11}(i) provides a first-order approximation of the curvature-related term in terms of the ratio of preconditioned gradient and momentum magnitudes. Along the discrete Adam iterates, this ratio and the coordinatewise quantity take the forms
\begin{equation*}
	\begin{aligned}
		\mu_k^{(i)}
		:=
		\frac{
			\|y_k^{(i)}\|
		}{
			\|x_k^{(i)}\|
		},\qquad
		\chi_k^{(i)}
		:=
		\frac{
			\|y_k\odot y_k^{(i)}\|^2
		}{
			\|y_k^{(i)}\|^2
		}.
	\end{aligned}
\end{equation*}
Combining the task-wise differences, we define a discrete indicator for $\omega_k$:
\begin{equation*}
	\label{eq:omega_hat}
	\begin{aligned}
		\hat{\omega}_k
		:=1+(1-\beta_1)
		\underbrace{
			\left(
			\mu_{k-1}^{(1)}
			-\mu_{k-1}^{(2)}
			\right)
		}_{\text{Term A}}+\frac{1-\beta_2}{2}
		\underbrace{
			\left(
			\chi_{k-1}^{(2)}
			-\chi_{k-1}^{(1)}
			\right)
		}_{\text{Term B}}.
	\end{aligned}
\end{equation*}
To assess the accuracy of $\hat{\omega}_k$ beyond the constant-step-size setting underlying its derivation, we revisit the Poisson benchmark with the exponentially decaying learning rate $\eta_k=10^{-3}\cdot 0.9^{k/1000}$. Figure~\ref{fig:omega_hat_comparison} compares $\omega_k$ and $\hat{\omega}_k$ over three randomly selected iteration intervals. Despite a systematic discrepancy in magnitude, $\hat{\omega}_k$ closely tracks the local trends and fluctuations of $\omega_k$ within \textbf{Phase I}.
Having observed that $\hat{\omega}_k$ captures the local variations of $\omega_k$ within Phase I, we next examine how task dominance in Adam's shared second-moment estimate relates to the two factors governing the step-to-step evolution of the task-contribution ratio $\rho_k$. Using the Poisson benchmark, we compare Adam without loss reweighting with Adam equipped with either Grad Norm or NTK-based adaptive loss reweighting. For this comparison, we introduce two auxiliary second-moment estimates, $v_k^{(1)}$ and $v_k^{(2)}$, obtained by accumulating the squared gradients of each task separately:
\begin{equation*}
	\begin{aligned}
		v^{(1)}_{k}
		= \beta_2 v^{(1)}_{k-1}
		+(1-\beta_2)(g_{k-1}^{(1)})^2, \qquad
		v^{(2)}_{k}
		= \beta_2 v^{(2)}_{k-1}
		+(1-\beta_2)(g_{k-1}^{(2)})^2,
	\end{aligned}
\end{equation*}
with $v_0^{(1)}=v_0^{(2)}=0$. Under the adaptive reweighted experiments, the task gradients in \eqref{eq:split} and the auxiliary recursions are replaced by their weighted counterparts:
\begin{equation*}
	\begin{aligned}
		g_k^{(1)}
		= \lambda_{bc}
		\nabla_\theta\mathcal{L}_{bc}(\theta_k), \qquad g_k^{(2)}
		= \lambda_{r}
		\nabla_\theta\mathcal{L}_{r}(\theta_k).
	\end{aligned}
\end{equation*}
The relative errors in Fig.~\ref{fig:secondmoment_comparison} reveal a pronounced dominance of the PDE-residual component under standard Adam. Specifically, the relative error between $v$ and $v^{(2)}$ remains near $10^{-5}$, whereas that between $v$ and $v^{(1)}$ stays close to unity. Both weighting schemes alleviate this dominance, as indicated by the increased discrepancy between $v$ and $v^{(2)}$. Under NTK weighting, this discrepancy increases in distinct steps following successive weight updates.

\begin{figure}[htbp]
	\centering
	\begin{subfigure}[b]{0.32\textwidth}
		\centering
		\includegraphics[width=\linewidth]{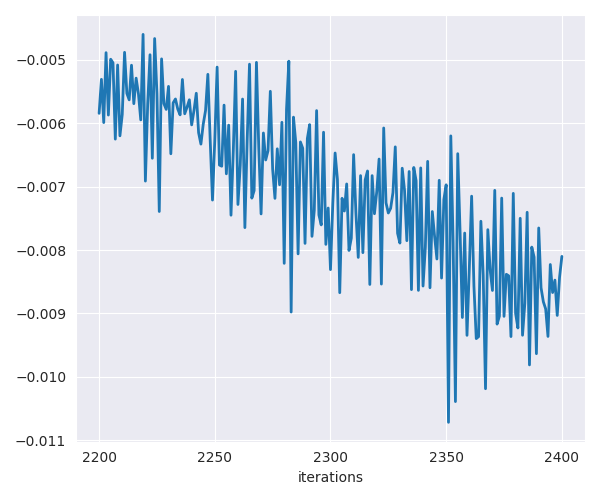}
		\caption{}
		\label{fig:term_a_adam}
	\end{subfigure}
	\hfill
	\begin{subfigure}[b]{0.32\textwidth}
		\centering
		\includegraphics[width=\linewidth]{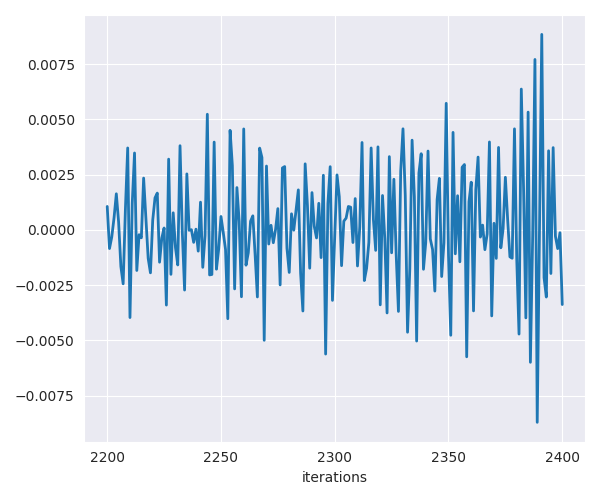}
		\caption{}
		\label{fig:term_a_gradnorm}
	\end{subfigure}
	\hfill
	\begin{subfigure}[b]{0.32\textwidth}
		\centering
		\includegraphics[width=\linewidth]{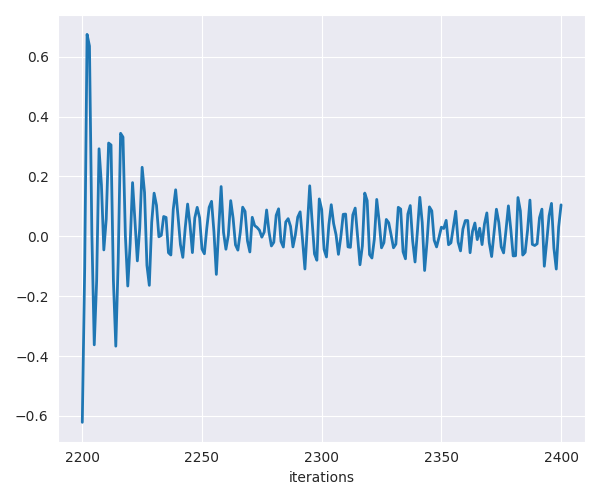}
		\caption{}
		\label{fig:term_a_ntk}
	\end{subfigure}
	\caption{
		Evolution of Term A,
		$\mu_k^{(1)}-\mu_k^{(2)}$,
		for Adam trained
		(a) without loss reweighting,
		(b) with gradient-norm-based loss reweighting,
		and (c) with NTK-based loss reweighting.
		The vertical scales differ across panels.
	}
	\label{fig:term_a_comparison}
\end{figure}

This change in the shared second-moment dynamics is accompanied by distinct behaviors of Term A, as illustrated in Fig.~\ref{fig:term_a_comparison}. Under standard Adam, $\mu_k^{(1)}-\mu_k^{(2)}$ remains negative and exhibits a progressively stronger negative bias. In contrast, both reweighting strategies induce repeated sign changes in this factor, replacing the persistent one-sided correction with alternating task-wise effects.

\begin{figure}[htbp]
	\centering
	\begin{subfigure}[b]{0.32\textwidth}
		\centering
		\includegraphics[width=\linewidth]{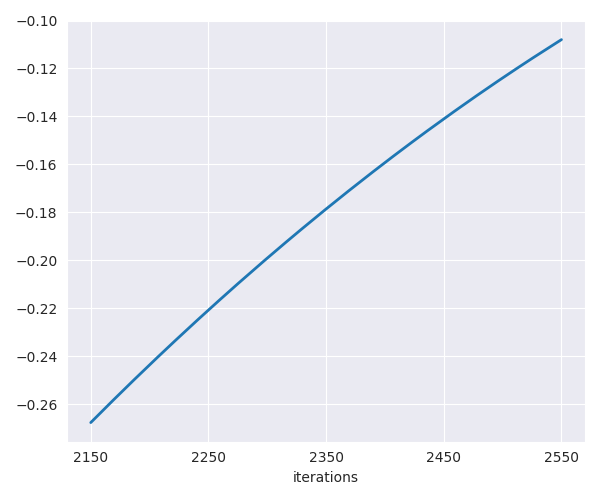}
		\caption{}
		\label{fig:term_b_adam}
	\end{subfigure}
	\hfill
	\begin{subfigure}[b]{0.32\textwidth}
		\centering
		\includegraphics[width=\linewidth]{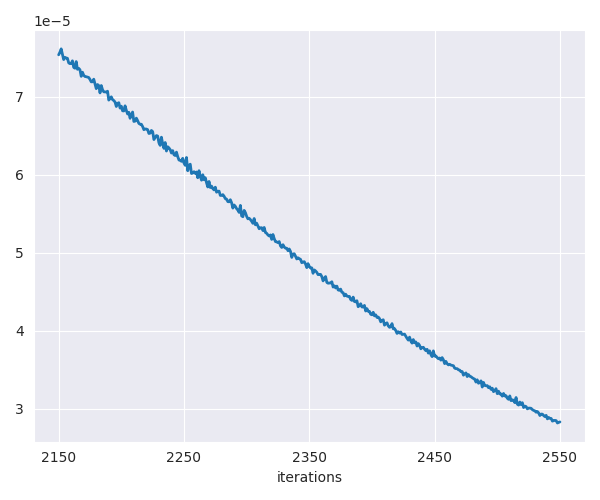}
		\caption{}
		\label{fig:term_b_gradnorm}
	\end{subfigure}
	\hfill
	\begin{subfigure}[b]{0.32\textwidth}
		\centering
		\includegraphics[width=\linewidth]{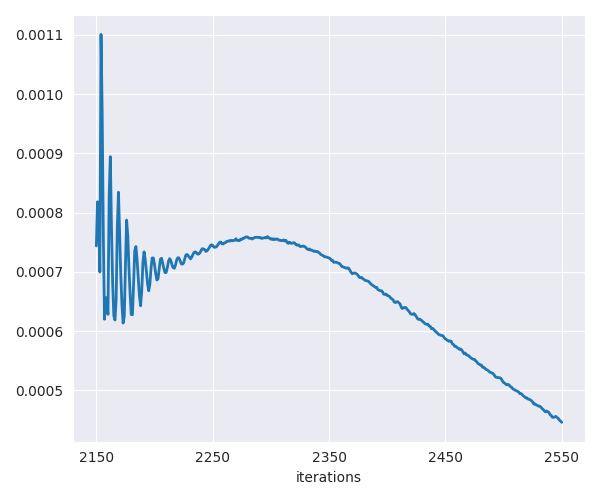}
		\caption{}
		\label{fig:term_b_ntk}
	\end{subfigure}
	\caption{
		Evolution of Term B,
		$\chi_k^{(2)}-\chi_k^{(1)}$,
		for Adam trained
		(a) without loss reweighting,
		(b) with gradient-norm-based loss reweighting,
		and (c) with NTK-based loss reweighting.
		The vertical scales differ across panels.
	}
	\label{fig:term_b_comparison}
\end{figure}

Term B exhibits a different response, as shown in Fig.~\ref{fig:term_b_comparison}. Under standard Adam, $\chi_k^{(2)}-\chi_k^{(1)}$ remains negative with a relatively large magnitude, reflecting a substantial weighted disparity in the coordinatewise energy distributions of the preconditioned task gradients. Both reweighting schemes reduce its absolute value by several orders of magnitude, substantially diminishing its influence on $\hat{\omega}_k$.

For $(\beta_1,\beta_2)=(0.9,0.999)$, Term A and Term B enter $\hat{\omega}_k$ with coefficients $0.1$ and $5\times10^{-4}$, respectively. Under standard Adam, the increasingly negative Term A becomes progressively dominant as Term B weakens, driving $\hat{\omega}_k$ further below unity, which suggests a sustained shift toward PDE-residual dominance. By contrast, the oscillatory behavior of Term A and the substantially diminished Term B under both weighting strategies suppress this persistent negative bias, favoring a more stable balance between the two tasks.

\subsection{Proposed Algorithm}
The preceding analysis highlights the influence of Adam's shared second-moment estimate on task-contribution dynamics. While loss-reweighting strategies modify this estimate by rescaling individual task gradients, AdaTask \cite{AdaTask_AAAI2023} instead maintains task-specific second-moment estimates:
\begin{equation}
	\label{eq:adatask}
	\begin{aligned}
		m_k^{(i)}
		&= \beta_1 m_{k-1}^{(i)}
		+(1-\beta_1)g_{k-1}^{(i)},\quad i=1,2,\\
		v_k^{(i)}
		&= \beta_2 v_{k-1}^{(i)}
		+(1-\beta_2)(g_{k-1}^{(i)})^2,
		\quad i=1,2,\\
		\hat m_k^{(i)}
		&= \frac{m_k^{(i)}}{1-\beta_1^k},
		\quad
		\hat v_k^{(i)}
		= \frac{v_k^{(i)}}{1-\beta_2^k}, \quad \; \; i=1,2,\\
		\theta_k
		&= \theta_{k-1}
		-\eta_k
		\sum_{i=1}^{2}
		\frac{\hat m_k^{(i)}}
		{\sqrt{\hat v_k^{(i)}}+\delta}.
	\end{aligned}
\end{equation}
Here, the two task-specific Adam updates are computed separately and then summed to update the shared parameters. MultiAdam \cite{pmlr-v202-yao23c}
employs the same task-wise Adam recursions but averages the resulting updates across loss groups, yielding
\begin{equation*}
	\theta_k
	=\theta_{k-1}
	-\frac{\eta_k}{2}
	\sum_{i=1}^{2}
	\frac{\hat m_k^{(i)}}
	{\sqrt{\hat v_k^{(i)}}+\delta}.
\end{equation*}
Hence, provided that the task grouping, initialization, and other Adam hyperparameters are identical, MultiAdam is equivalent to AdaTask with a rescaled learning rate. We therefore use AdaTask as the representative of both methods in the subsequent numerical experiments.

Although AdaTask employs task-specific second-moment estimates, its update admits an equivalent representation that retains Adam's shared second-moment estimate by appropriately rescaling the task-specific momenta. Specifically, with $M_0^{(i)}=\hat v_0^{(i)}=\hat v_0=0$ by convention, the modulated momentum of each task satisfies the following recursion for $k\geq1$:
\begin{equation*}
	\begin{aligned}
		M_k^{(i)}
		={}
		\beta_1
		{
			\frac{\sqrt{\hat v_{k-1}^{(i)}}+\delta}
			{\sqrt{\hat v_k^{(i)}}+\delta}
		}
		{
			\frac{\sqrt{\hat v_k}+\delta}
			{\sqrt{\hat v_{k-1}}+\delta}
		}
		M_{k-1}^{(i)}
		+(1-\beta_1)
		(\frac{\sqrt{\hat v_k}+\delta}
		{\sqrt{\hat v_k^{(i)}}+\delta}
		g_{k-1}^{(i)}),\quad \;  i=1,2.
	\end{aligned}
\end{equation*}
Replacing the task-specific momenta in the bias correction and parameter update of \eqref{eq:split} with $M_k^{(i)}$ recovers exactly the AdaTask update in \eqref{eq:adatask}. In this formulation, task-dependent rescaling acts on both the accumulated momentum and the incoming gradients. Rather than rescaling the accumulated momentum at every iteration, we consider applying task-dependent modulation directly to the incoming gradients. This choice is motivated by the fact that task gradients drive the evolution of Adam's momentum, while its second-moment estimate may respond sluggishly to instantaneous squared gradients in the decoupling regime identified by Bai et al.~\cite{pmlr-v306-bai26o}. Accordingly, we propose Gradient-Modulated Adam (GMAdam), defined by the following iteration:
\begin{equation}
	\label{eq:proposed_adatask}
	\begin{aligned}
		v_k^{(i)}
		&= \beta_2 v_{k-1}^{(i)}
		+(1-\beta_2)(g_{k-1}^{(i)})^2,
		\quad i=1,2,\\
		v_k
		&= \beta_2 v_{k-1}
		+(1-\beta_2)
		\left(
		\sum_{i=1}^{2}g_{k-1}^{(i)}
		\right)^2,\\
		m_k
		&= \beta_1m_{k-1}
		+(1-\beta_1)
		\sum_{i=1}^{2}
		\sqrt{
			\frac{v_k}{v_k^{(i)}+\delta_1}
		}
		\odot g_{k-1}^{(i)},\\
		\hat m_k
		&= \frac{m_k}{1-\beta_1^k},
		\qquad
		\hat v_k
		= \frac{v_k}{1-\beta_2^k},\\
		\theta_k
		&= \theta_{k-1}
		-\eta_k
		\frac{\hat m_k}
		{\sqrt{\hat v_k}+\delta}.
	\end{aligned}
\end{equation}
The positive parameter $\delta_1$ prevents excessively small task-specific second-moment denominators, thereby controlling the magnitude of the gradient modulation factors. In the numerical experiments, $\delta_1$ takes values of $10^{-8}$ or $10^{-12}$, with the specific setting reported for each test case, while $\delta$ is fixed at $10^{-8}$ throughout.

From a computational perspective, the proposed method requires $(n+2)d$ optimizer-state scalars for $n$ tasks and $d$ trainable parameters, compared with $2nd$ for AdaTask. This yields a reduction of $(n-2)d$, which grows linearly with the number of tasks. We evaluate GMAdam against AdaTask and two loss-reweighting baselines, namely Grad Norm and NTK-based weighting, on four PDE benchmarks: the Poisson, wave, and Allen--Cahn equations, and lid-driven cavity flow ($Re=100,400$). The two weighting schemes follow the updated formulations and publicly available implementation of Wang et al.~\cite{wang2023expertsguidetrainingphysicsinformed}. Across all test configurations, GMAdam achieves training times within $3 \% $ of those of AdaTask, while consistently yielding lower relative $L^2$ errors. Moreover, it attains the lowest errors among all four methods for the Poisson, wave, and Allen--Cahn equations.

\section{Numerical Results}
All numerical experiments were conducted on an NVIDIA GeForce RTX 5070 Laptop GPU using single-precision (float32) arithmetic. Throughout the experiments, the decay parameters were fixed at $(\beta_1,\beta_2)=(0.9,0.999)$ for all considered methods.

\subsection{1D Poisson Equation}
\label{subsec:poisson_results}
We revisit the high-frequency Poisson problem considered in Section~3.1,
\begin{equation*}
	\label{eq:poisson_numerical}
	\begin{aligned}
		-\Delta u(x)
		&= 4\pi^2\sin(2\pi x)
		+10^4\pi^2\sin(100\pi x),
		\quad x\in(0,1),\\
		u(0)&=u(1)=0,
	\end{aligned}
\end{equation*}
whose exact solution is $u(x)=\sin(2\pi x)+\sin(100\pi x)$. We approximate $u$ by the four-layer Fourier-feature MLP employed in Section~3.1, retaining the same embedding scale $\sigma=13.5$. For AdaTask and GMAdam, the training objective comprises the boundary-condition and PDE-residual losses directly,
\begin{equation*}
	\begin{aligned}
		\mathcal{L}(\theta)
		&=\mathcal{L}_{bc}(\theta)+\mathcal{L}_{r}(\theta)\\
		&= \frac{1}{N_{bc}}\sum_{i=1}^{N_{bc}} |u_\theta(x^i_{bc})|^2 +\frac{1}{N_r}\sum_{j=1}^{N_r}
		\left|
		-\Delta u_\theta(x_r^j)
		-f(x_r^j)
		\right|^2,
	\end{aligned}
\end{equation*}
where $f(x)=4\pi^2\sin(2\pi x)+10^4\pi^2\sin(100\pi x)$. $\mathcal{L}_{bc}$ and $\mathcal{L}_{r}$ are treated as two distinct tasks in their corresponding updates. By contrast, the Grad Norm and NTK-based reweighting schemes employ Adam with the iteration-dependent objective:
\begin{equation*}
	\mathcal{L}(\theta)
	=\lambda_{bc}\mathcal{L}_{bc}(\theta)
	+\lambda_{r}\mathcal{L}_{r}(\theta).
\end{equation*}
As recommended by Wang et al.~\cite{wang2023expertsguidetrainingphysicsinformed}, both reweighting schemes use an EMA coefficient of $\alpha=0.9$ and a weight-update interval of $T=1000$ iterations. These settings are maintained throughout the subsequent experiments. For each endpoint, 128 boundary samples are used, giving $N_{bc}=256$ in total, while $N_r=4{,}000$ interior collocation points are generated using Latin hypercube sampling (LHS)~\cite{Stein01051987} and resampled every 20 iterations. All four methods are trained for $10^5$ iterations with an initial learning rate of $5\times10^{-4}$ and exponential learning-rate decay with a decay-step parameter of 2,000 iterations. GMAdam uses $\delta_1=10^{-8}$. Prediction accuracy is evaluated using the relative $L^2$ error computed over $N_{\mathrm{test}}=10{,}000$ uniformly spaced test points in $[0,1]$.

\begin{figure}[htbp]
	\centering
	\includegraphics[width=\textwidth]{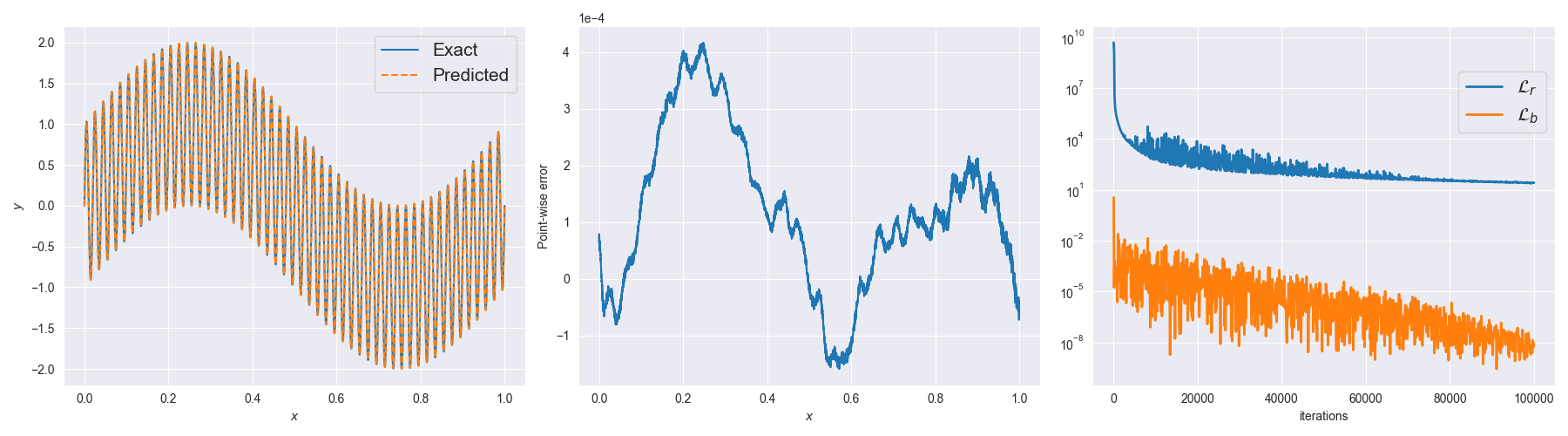}
	\caption{One-dimensional Poisson equation: results obtained with
		GMAdam after $10^5$ iterations.
		Left: exact and predicted solutions.
		Middle: signed pointwise error $u-u_\theta$.
		Right: evolution of the boundary loss $\mathcal{L}_{bc}$
		and PDE-residual loss $\mathcal{L}_r$ (logarithmic scale).}
	\label{fig:poisson_solution}
\end{figure}

\begin{figure}[H]
	\centering
	\begin{subfigure}[t]{0.245\textwidth}
		\centering
		\includegraphics[width=\linewidth]
		{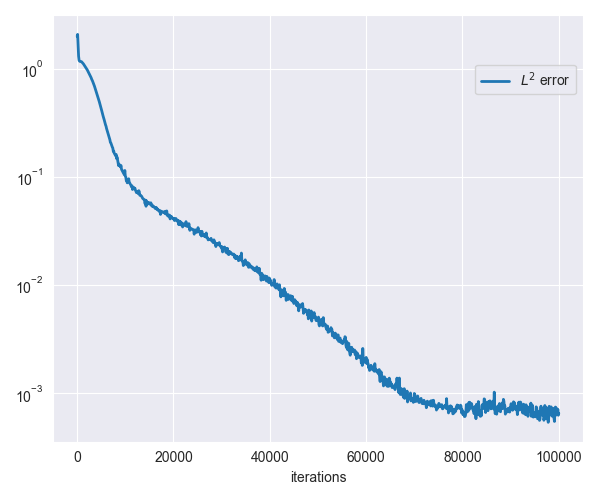}
		\caption{}
		\label{fig:poisson_l2_adatask}
	\end{subfigure}
	\begin{subfigure}[t]{0.245\textwidth}
		\centering
		\includegraphics[width=\linewidth]
		{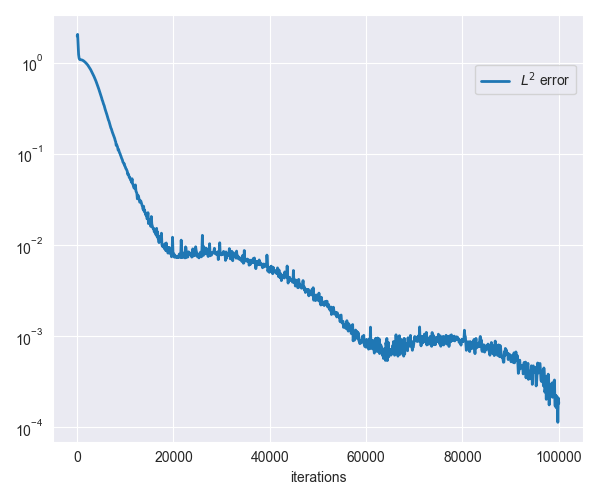}
		\caption{}
		\label{fig:poisson_l2_ours}
	\end{subfigure}
	\begin{subfigure}[t]{0.245\textwidth}
		\centering
		\includegraphics[width=\linewidth]
		{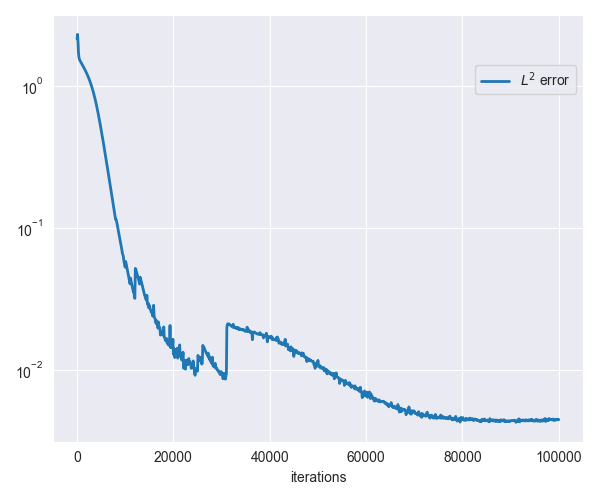}
		\caption{}
		\label{fig:poisson_l2_gradnorm}
	\end{subfigure}
	\begin{subfigure}[t]{0.245\textwidth}
		\centering
		\includegraphics[width=\linewidth]
		{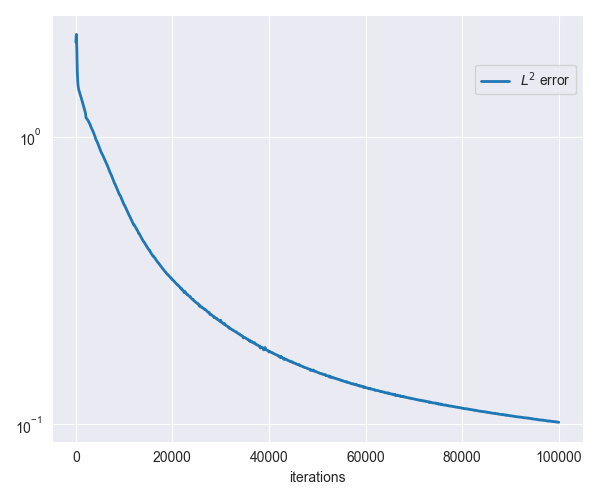}
		\caption{}
		\label{fig:poisson_l2_ntk}
	\end{subfigure}
	\caption{
		Evolution of the relative $L^2$ error
		during training for
		(a) AdaTask,
		(b) GMAdam,
		(c) Adam with Grad Norm
		loss weighting, and
		(d) Adam with NTK-based loss weighting.
		The vertical axes are logarithmic.
	}
	\label{fig:poisson_l2_histories}
\end{figure}

Figure~\ref{fig:poisson_solution} presents the numerical solution obtained with the proposed method, together with its pointwise error and the histories of the two unweighted loss components. The predicted solution closely reproduces both frequency components, with pointwise discrepancies of order $10^{-4}$.
The quantitative comparison in Table~\ref{tab:poisson_comparison} shows that GMAdam achieves the lowest relative $L^2$ error among the four approaches. In particular, it reduces the error of AdaTask by approximately $71.9\%$, with a $2.8\%$ increase in training time.
\begin{table}[H]
	\centering
	\caption{Relative $L^2$ errors and training times
		for the 1D Poisson equation after $10^5$ iterations.}
	\label{tab:poisson_comparison}
	\begin{tabular}{lcc}
		\toprule
		Method & Relative $L^2$ error & training time (s) \\
		\midrule
		AdaTask                 & $6.56\times10^{-4}$ & 340.49 \\
		GMAdam         & $\mathbf{1.84\times10^{-4}}$ & 349.89 \\
		Grad Norm & $4.47\times10^{-3}$ & 315.15 \\
		NTK           & $1.02\times10^{-1}$ & 550.96 \\
		\bottomrule
	\end{tabular}
\end{table}
The error histories of the four methods are compared in Fig.~\ref{fig:poisson_l2_histories}.
\subsection{1D Wave Equation}
We next assess GMAdam on the following initial-boundary value problem for the one-dimensional wave equation:
\begin{equation*}
	\begin{aligned}
		&u_{tt}(x,t)-25u_{xx}(x,t)=0,
		\qquad (x,t)\in(0,1)\times(0,1),\\
		&u(0,t)=u(1,t)=0,
		\qquad t\in[0,1],\\
		&u(x,0)=\sin(\pi x)+\sin(4\pi x),
		\qquad x\in[0,1],\\
		&u_t(x,0)=0,
		\qquad x\in[0,1].
	\end{aligned}
\end{equation*}
The prescribed initial displacement excites two standing-wave modes, yielding the exact solution
\begin{equation*}
	u(x,t)=\sin(\pi x)\cos(5\pi t)
	+\sin(4\pi x)\cos(20\pi t).
\end{equation*}
We employ the modified MLP architecture of Wang et al.~\cite{wang2021understanding}, with one gated hidden layer of width 256 and $\tanh$ activation. To accommodate the distinct spectral characteristics of the spatial and temporal variables, the network incorporates an anisotropic Fourier feature embedding (AFFE) following Jabareen et al.~\cite{10.1007/978-3-032-06774-6_1}. Compared with the single-scale and multi-scale Fourier feature architectures of Wang et al.~\cite{wang2021eigenvector}, the embedding assigns independent frequency scales to individual input coordinates within a single branch. Specifically, for an input $\boldsymbol{x}\in\mathbb{R}^{p_0}$, the mapping is defined by
\begin{equation}
	\label{eq:anisotropic_fourier_embedding}
	\gamma(\boldsymbol{x})=
	\begin{bmatrix}
		\sin\left(BD_{\sigma}\boldsymbol{x}\right)\\
		\cos\left(BD_{\sigma}\boldsymbol{x}\right)
	\end{bmatrix},
\end{equation}
where $B\in\mathbb{R}^{\frac{p_1}{2}\times p_0}$ is a fixed random matrix with independent standard Gaussian entries, and $D_{\sigma}=\operatorname{diag}(\sigma_1,\ldots,\sigma_{p_0})$ specifies the Fourier scale associated with each input coordinate. Here, the spatial and temporal Fourier scales are set to $\sigma_{\mathbf{x}}=2$ and $\sigma_t=6$, respectively. 

For AdaTask and GMAdam, the training objective consists of four loss components,
\begin{equation*}
	\mathcal{L}(\theta)
	=\mathcal{L}_{bc,1}(\theta)
	+\mathcal{L}_{bc,2}(\theta)
	+\mathcal{L}_{ic}(\theta)
	+\mathcal{L}_{r}(\theta),
\end{equation*}
where
\begin{equation*}
	\begin{aligned}
		\mathcal{L}_{bc,1}(\theta)
		={}&
		\frac{1}{N_{bc}}\sum_{i=1}^{N_{bc}}
		|u_\theta(0,t_l^i)|^2
		+\frac{1}{N_{bc}}\sum_{i=1}^{N_{bc}}
		|u_\theta(1,t_u^i)|^2,\\
		\mathcal{L}_{bc,2}(\theta)
		={}&
		\frac{1}{N_{ic}}\sum_{i=1}^{N_{ic}}
		|u_\theta(x_{ic}^i,0)-u_0(x_{ic}^i)|^2,\\
		\mathcal{L}_{ic}(\theta)
		={}&
		\frac{1}{N_{ic}}\sum_{i=1}^{N_{ic}}
		\left|
		\frac{\partial u_{\theta}}{\partial t}(x_{ic}^i,0)
		\right|^2,\\
		\mathcal{L}_{r}(\theta)
		={}&
		\frac{1}{N_r}\sum_{i=1}^{N_r}
		\left|
		\frac{\partial^2 u_{\theta}}{\partial t^2}(x_r^i,t_r^i)
		-25\frac{\partial^2 u_{\theta}}{\partial x^2}(x_r^i,t_r^i)
		\right|^2,
	\end{aligned}
\end{equation*}
with $u_0(x)=\sin(\pi x)+\sin(4\pi x)$. For Grad Norm and NTK-based reweighting schemes, the boundary-condition and initial-displacement losses are combined into a single weighted group, leading to the following objective optimized by Adam:
\begin{equation*}
	\mathcal{L}(\theta)
	=\lambda_{bc}
	\left(
	\mathcal{L}_{bc,1}(\theta)
	+\mathcal{L}_{bc,2}(\theta)
	\right)
	+\lambda_{ic}\mathcal{L}_{ic}(\theta)
	+\lambda_r\mathcal{L}_r(\theta).
\end{equation*}
At each training iteration, LHS is used to sample $N_{bc}=32$ points on each spatial boundary, $N_{ic}=64$ points at $t=0$, and $N_r=2{,}000$ interior collocation points. Training is conducted for $8\times10^4$ iterations using an exponentially decaying learning rate $\eta_k=2\times10^{-4}(0.9)^{k/1000}$. For GMAdam, $\delta_1$ is fixed at $10^{-12}$. The relative $L^2$ error is measured on a uniform $500\times500$ grid covering the space-time domain $[0,1]^2$.

\begin{figure}[htbp]
	\centering
	\includegraphics[width=\textwidth]{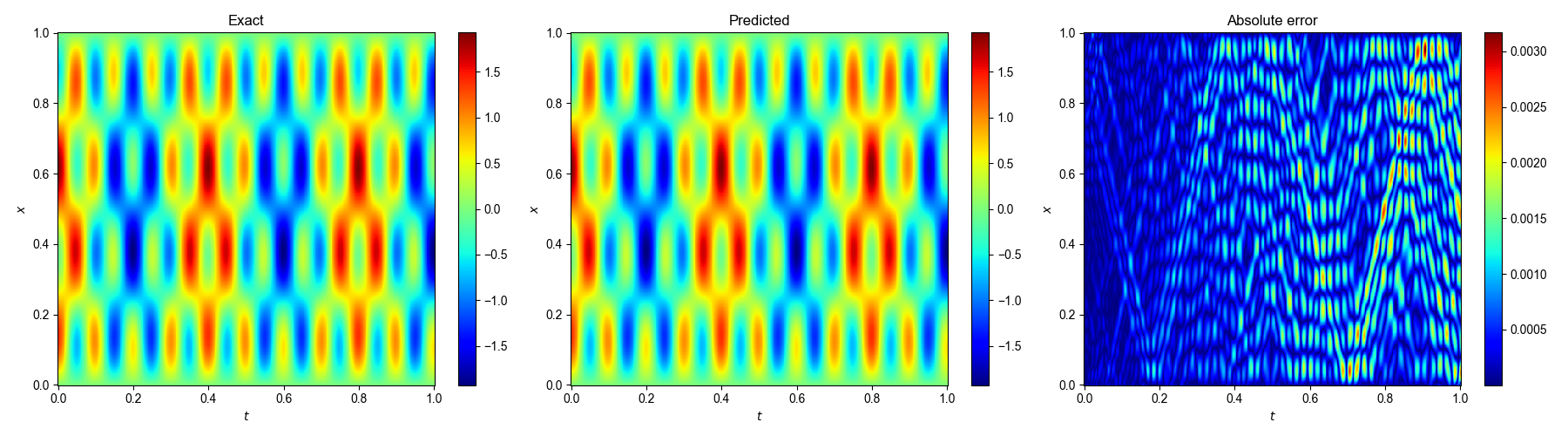}
	\caption{One-dimensional wave equation: results obtained with GMAdam after $8\times10^4$ iterations. Left: exact solution. Middle: predicted solution. Right: pointwise absolute error $|u-u_\theta|$.}
	\label{fig:wave_solution}
\end{figure}

The numerical solution obtained with GMAdam is compared with the exact solution in Fig.~\ref{fig:wave_solution}. The predicted field accurately captures the oscillatory structure of both standing-wave modes, with pointwise absolute errors on the order of $10^{-3}$ or smaller over most of the space-time domain.

Table~\ref{tab:wave_comparison} summarizes the relative $L^2$ errors and computational times of the four methods. Although all approaches achieve comparable accuracy, the proposed method attains the lowest error of $9.86\times10^{-4}$, representing a $19.7\%$ reduction relative to AdaTask while requiring approximately $1.8\%$ less training time. The Grad Norm scheme incurs the lowest computational cost, whereas NTK-based weighting requires the longest training time.

\begin{table}[H]
	\centering
	\caption{Relative $L^2$ errors and training times for the one-dimensional wave equation after $8\times10^4$ iterations.}
	\label{tab:wave_comparison}
	\begin{tabular}{lcc}
		\toprule
		Method & Relative $L^2$ error & Training time (s) \\
		\midrule
		AdaTask          & $1.2273\times10^{-3}$ & 1249.63 \\
		GMAdam  & $\mathbf{9.8599\times10^{-4}}$ & 1227.44 \\
		Grad Norm        & $1.0139\times10^{-3}$ & 1135.06 \\
		NTK              & $1.0544\times10^{-3}$ & 1502.02 \\
		\bottomrule
	\end{tabular}
\end{table}

\begin{figure}[H]
	\centering
	\begin{subfigure}[t]{0.245\textwidth}
		\centering
		\includegraphics[width=\linewidth]{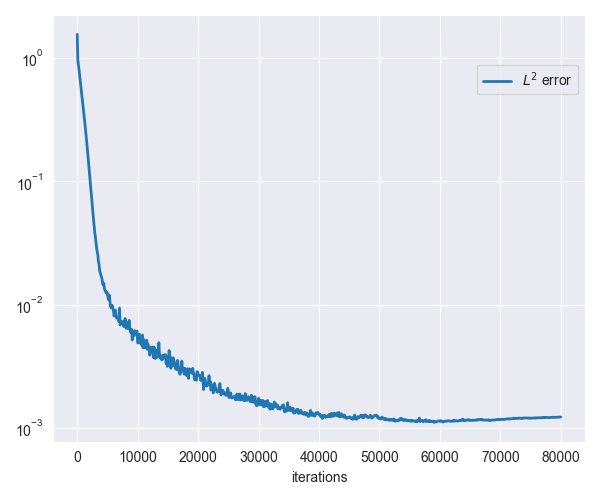}
		\caption{}
		\label{fig:wave_l2_adatask}
	\end{subfigure}
	\begin{subfigure}[t]{0.245\textwidth}
		\centering
		\includegraphics[width=\linewidth]{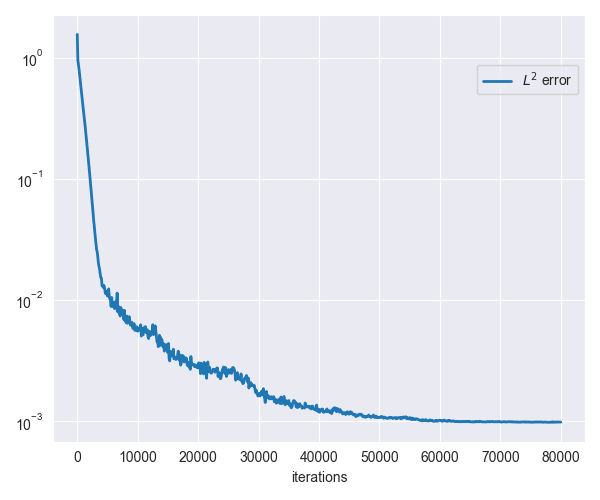}
		\caption{}
		\label{fig:wave_l2_ours}
	\end{subfigure}
	\begin{subfigure}[t]{0.245\textwidth}
		\centering
		\includegraphics[width=\linewidth]{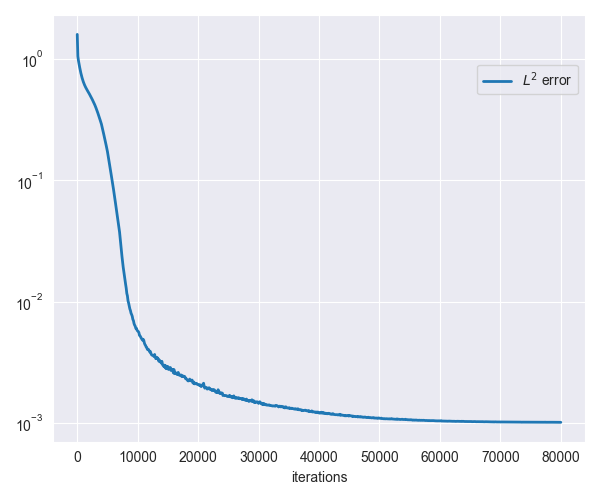}
		\caption{}
		\label{fig:wave_l2_gradnorm}
	\end{subfigure}
	\begin{subfigure}[t]{0.245\textwidth}
		\centering
		\includegraphics[width=\linewidth]{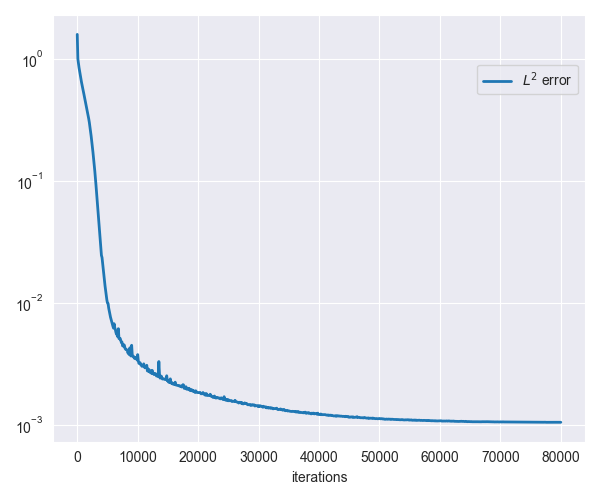}
		\caption{}
		\label{fig:wave_l2_ntk}
	\end{subfigure}
	\caption{Evolution of the relative $L^2$ error during training for (a) AdaTask, (b) GMAdam, (c) Adam with Grad Norm, and (d) Adam with NTK-based loss weighting. The vertical axes are logarithmic.}
	\label{fig:wave_l2_histories}
\end{figure}

The corresponding convergence histories are displayed in Fig.~\ref{fig:wave_l2_histories}. All four methods exhibit a rapid initial reduction in relative $L^2$ error, followed by progressively slower convergence. AdaTask shows a slight increase in error toward the end of training, whereas  no comparable increase is observed for GMAdam. Both reweighting schemes exhibit relatively smooth late-stage trajectories, approaching final accuracy levels close to that of the proposed method.

\subsection{Allen-Cahn Equation}
As a nonlinear benchmark, we consider the one-dimensional Allen--Cahn equation subject to periodic boundary conditions:
\begin{equation*}
	\begin{aligned}
		&u_t-\nu u_{xx}+f(u)=0,
		\qquad (x,t)\in(-1,1)\times(0,1),\\
		&u(x,0)=x^2\cos(\pi x),
		\qquad x\in[-1,1],\\
		&u(-1,t)=u(1,t),\quad
		u_x(-1,t)=u_x(1,t),
		\qquad t\in[0,1],
	\end{aligned}
\end{equation*}
where $\nu=10^{-4}$ and $f(u)=5(u^3-u)$. The network comprises two hidden layers of width 256 with $\tanh$ activation, preceded by an AFFE mapping with temporal and spatial scales $\sigma_t=2$ and $\sigma_x=10$, respectively.

For AdaTask and GMAdam, the loss function is decomposed into four components:
\begin{equation*}
	\mathcal{L}(\theta)
	=\mathcal{L}_{bc,1}(\theta)
	+\mathcal{L}_{bc,2}(\theta)
	+\mathcal{L}_{ic}(\theta)
	+\mathcal{L}_{r}(\theta),
\end{equation*}
where
\begin{equation*}
	\begin{aligned}
		\mathcal{L}_{bc,1}(\theta)
		={}&
		\frac{1}{N_{bc}}\sum_{i=1}^{N_{bc}}
		|u_\theta(-1,t_{bc}^i)
		-u_\theta(1,t_{bc}^i)|^2,\\
		\mathcal{L}_{bc,2}(\theta)
		={}&
		\frac{1}{N_{bc}}\sum_{i=1}^{N_{bc}}
		\left|
		\frac{\partial u_\theta}{\partial x}(-1,t_{bc}^i)
		-\frac{\partial u_\theta}{\partial x}(1,t_{bc}^i)
		\right|^2,\\
		\mathcal{L}_{ic}(\theta)
		={}&
		\frac{1}{N_{ic}}\sum_{i=1}^{N_{ic}}
		|u_\theta(x_{ic}^i,0)-u_0(x_{ic}^i)|^2,\\
		\mathcal{L}_{r}(\theta)
		={}&
		\frac{1}{N_r}\sum_{i=1}^{N_r}
		\left|
		\frac{\partial u_\theta}{\partial t}(x_r^i,t_r^i)
		-\nu\frac{\partial^2 u_\theta}{\partial x^2}(x_r^i,t_r^i)
		+f\bigl(u_\theta(x_r^i,t_r^i)\bigr)
		\right|^2,
	\end{aligned}
\end{equation*}
with $u_0(x)=x^2\cos(\pi x)$. Grad Norm and NTK-based schemes instead employ a two-group loss formulation:
\begin{equation*}
	\mathcal{L}(\theta)
	=\lambda_{bc}
	\left(
	\mathcal{L}_{bc,1}(\theta)
	+\mathcal{L}_{bc,2}(\theta)
	+\mathcal{L}_{ic}(\theta)
	\right)
	+\lambda_r\mathcal{L}_r(\theta).
\end{equation*}
At every training iteration, $N_{bc}=50$ temporal points for enforcing both periodic boundary conditions, $N_{ic}=300$ initial-condition points, and $N_r=4{,}000$ interior collocation points are sampled uniformly at random from their respective domains. Training proceeds for $1.5 \times 10^5$ iterations with an exponentially decaying learning rate $\eta_k=10^{-3}(0.9)^{k/2000}$. The gradient-modulation parameter in GMAdam is set to $\delta_1=10^{-12}$. The relative $L^2$ error is evaluated against the numerical reference solution of Raissi et al.~\cite{raissi2019physics}.

\begin{figure}[htbp]
	\centering
	\includegraphics[width=\textwidth]{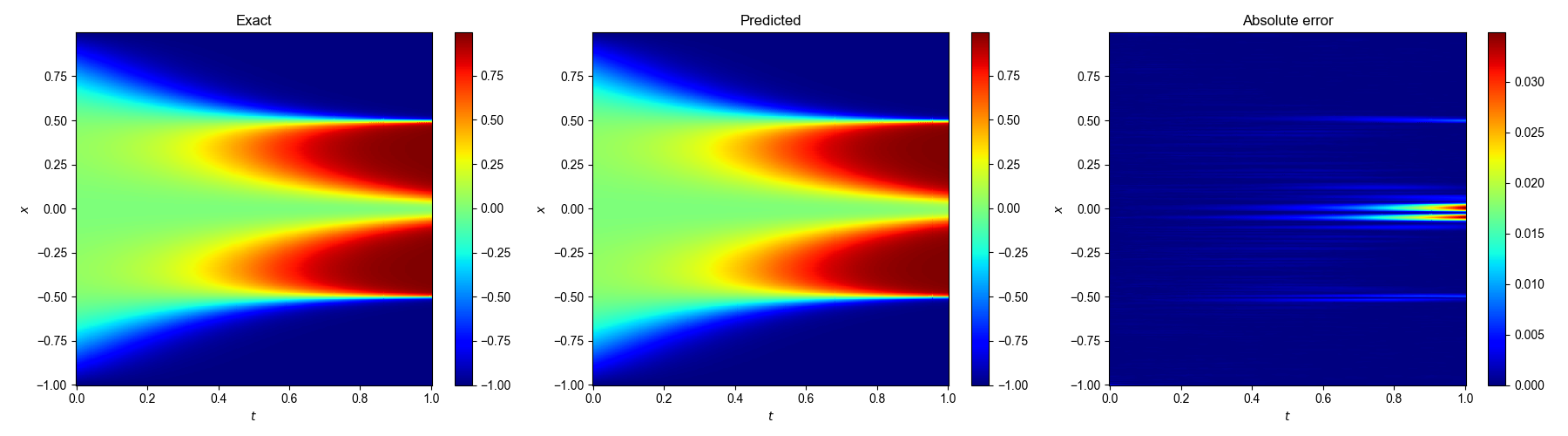}
	\caption{One-dimensional Allen--Cahn equation: reference solution (left), prediction obtained with GMAdam (middle), and pointwise absolute error $|u-u_\theta|$ (right) after $1.5\times10^5$ iterations.}
	\label{fig:ac_solution}
\end{figure}

Figure~\ref{fig:ac_solution} presents the numerical prediction obtained with the proposed method, together with the reference solution and the corresponding pointwise absolute error. The predicted field closely reproduces the reference solution, with the largest pointwise errors observed near $(x,t)=(0,1)$.

The quantitative results reported in Table~\ref{tab:ac_comparison} indicate that GMAdam yields the lowest relative $L^2$ error of $2.85\times10^{-3}$, improving upon AdaTask by approximately $21.0\%$ while reducing training time by $1.0\%$. The two reweighting schemes achieve comparable accuracy, with Grad Norm weighting requiring the shortest training time and NTK-based weighting incurring the highest computational cost.

\begin{table}[H]
	\centering
	\caption{Relative $L^2$ errors and training times for the one-dimensional Allen--Cahn equation after $1.5\times10^5$ iterations.}
	\label{tab:ac_comparison}
	\begin{tabular}{lcc}
		\toprule
		Method & Relative $L^2$ error & Training time (s) \\
		\midrule
		AdaTask         & $3.6097\times10^{-3}$ & 1626.82 \\
		GMAdam & $\mathbf{2.8518\times10^{-3}}$ & 1610.42 \\
		Grad Norm       & $2.9500\times10^{-3}$ & 1487.00 \\
		NTK             & $3.2818\times10^{-3}$ & 2164.36 \\
		\bottomrule
	\end{tabular}
\end{table}

\begin{figure}[H]
	\centering
	\begin{subfigure}[t]{0.24\textwidth}
		\centering
		\includegraphics[width=\linewidth]{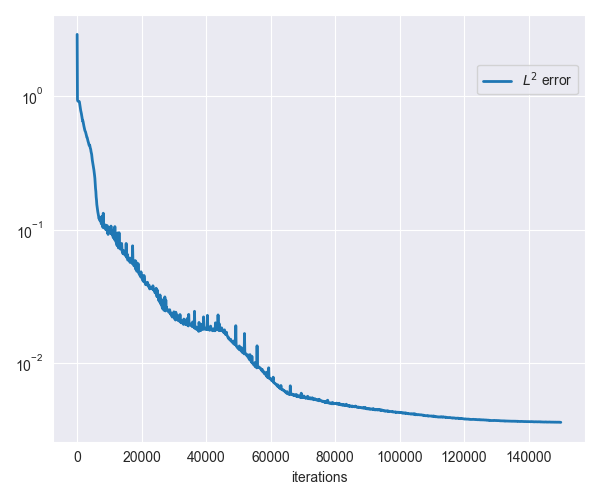}
		\caption{}
		\label{fig:ac_l2_adatask}
	\end{subfigure}
	\hfill
	\begin{subfigure}[t]{0.24\textwidth}
		\centering
		\includegraphics[width=\linewidth]{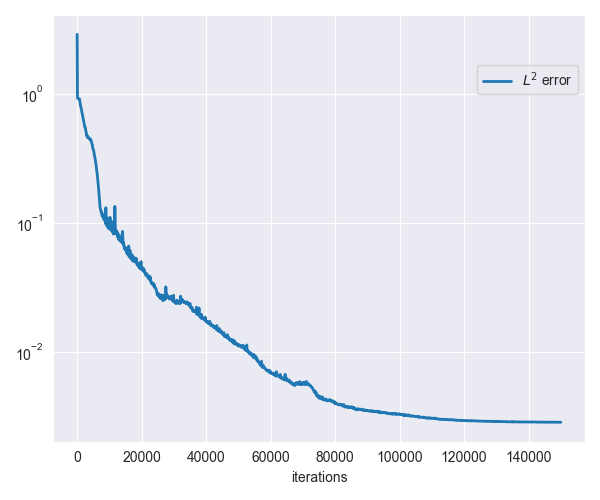}
		\caption{}
		\label{fig:ac_l2_ours}
	\end{subfigure}
	\hfill
	\begin{subfigure}[t]{0.24\textwidth}
		\centering
		\includegraphics[width=\linewidth]{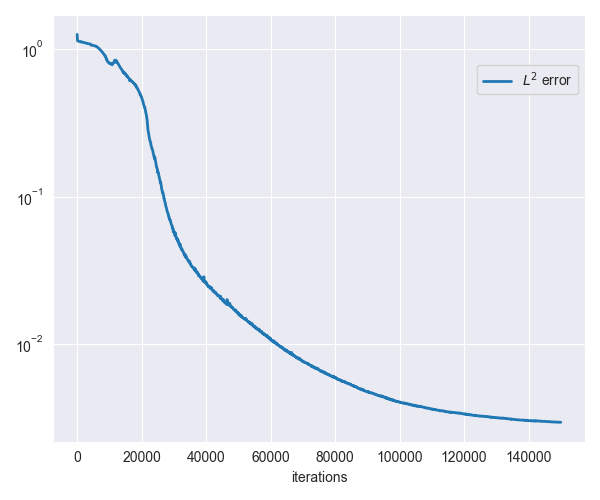}
		\caption{}
		\label{fig:ac_l2_gradnorm}
	\end{subfigure}
	\hfill
	\begin{subfigure}[t]{0.24\textwidth}
		\centering
		\includegraphics[width=\linewidth]{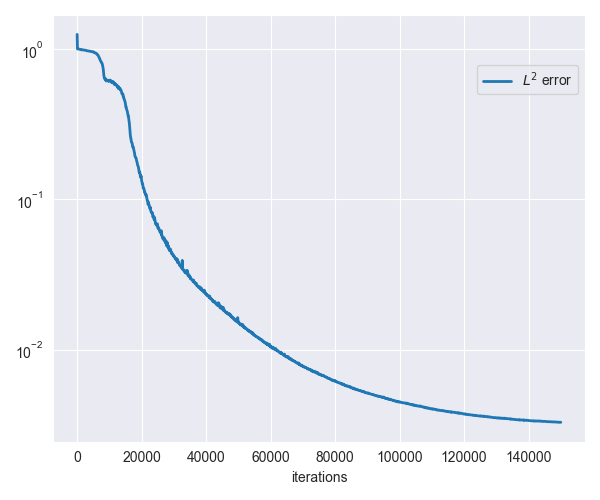}
		\caption{}
		\label{fig:ac_l2_ntk}
	\end{subfigure}
	\caption{Evolution of the relative $L^2$ error over $1.5\times10^5$ training iterations for (a) AdaTask, (b) GMAdam, (c) Adam with Grad Norm, and (d) Adam with NTK-based loss weighting. The vertical axes are logarithmic.}
	\label{fig:ac_l2_histories}
\end{figure}

Figure~\ref{fig:ac_l2_histories} compares the evolution of the relative $L^2$ error across the four methods. While all approaches exhibit an overall reduction in error, their transient behaviors differ. AdaTask displays noticeable fluctuations during the intermediate stages of training, which are less pronounced for the proposed method. In comparison, the Grad Norm and NTK-based schemes exhibit a more gradual initial decrease, followed by relatively smooth error reduction. The four trajectories eventually reach the same order of accuracy, with the proposed method attaining the smallest final error.

\subsection{Lid-Driven Cavity Flow}
The two-dimensional steady lid-driven cavity flow serves as the final benchmark, governed by the dimensionless incompressible Navier--Stokes equations:
\begin{equation*}
	\label{eq:ldc_ns}
	\begin{aligned}
		(\boldsymbol{u}\cdot\nabla)\boldsymbol{u}
		+\nabla p
		-\frac{1}{Re}\Delta\boldsymbol{u}
		&=\boldsymbol{0}, \qquad \text{in } \Omega, \\
		\nabla\cdot\boldsymbol{u} 
		&=0, \qquad   \text{in } \Omega.
 	\end{aligned}
\end{equation*}
Here, $\Omega=(0,1)^2$ specifies the unit-square domain, and $\boldsymbol{u}=(u,v)^{\mathsf T}$, $p$, and $Re$ denote the velocity field, pressure, and Reynolds number, respectively. The upper lid moves horizontally at unit speed, while the remaining walls are stationary, subject to the boundary conditions.
\begin{equation*}
	\label{eq:ldc_bc}
	\begin{aligned}
		\boldsymbol{u}(x,1)
		&=(1,0)^{\mathsf T},\\
		\boldsymbol{u}(x,0)
		&=\boldsymbol{0},\\
		\boldsymbol{u}(0,y)
		&=\boldsymbol{u}(1,y)=\boldsymbol{0}.
	\end{aligned}
\end{equation*}
This benchmark is investigated at two Reynolds numbers, $Re=100$ and $Re=400$.

The $Re=100$ configuration uses a modified MLP with one gated hidden layer of width 256 and $\tanh$ activation to jointly approximate the velocity and pressure fields. The network inputs are transformed using a single-scale Fourier feature embedding with $\sigma=2$. AdaTask and GMAdam treat the following two loss components as separate tasks:
\begin{equation*}
	\begin{aligned}
		\mathcal{L}_{r}(\theta)
		&=
		\mathcal{L}_{u}(\theta)
		+\mathcal{L}_{v}(\theta)
		+\mathcal{L}_{c}(\theta),\\
		\mathcal{L}_{bc}(\theta)
		&=
		\mathcal{L}_{bc}^{u}(\theta)
		+\mathcal{L}_{bc}^{v}(\theta).
	\end{aligned}
\end{equation*}
The individual contributions to $\mathcal{L}_r$ are evaluated as
\begin{equation*}
	\mathcal{L}_{q}(\theta)
	=
	\frac{1}{N_r}
	\sum_{i=1}^{N_r}
	\left|r_{q,\theta}(x_r^i,y_r^i)\right|^2,
	\quad \; \; q\in\{u,v,c\},
\end{equation*}
where $q \in  \{u,v\} $ indexes the two momentum equations, while $q=c$ refers to the incompressibility constraint. A fresh set of $N_r=2{,}000$ interior collocation points is drawn uniformly at random from $\Omega$ for each training iteration. The remaining component, $\mathcal{L}_{bc}$, penalizes deviations from the prescribed velocities on the cavity walls. Its two constituent terms are given by
\begin{equation*}
	\begin{aligned}
		\mathcal{L}_{bc}^{u}(\theta)
		&=\frac{1}{N_{bc}}\sum_{i=1}^{N_{bc}}
		\Big[
		|u_\theta(x_t^i,1)-u_{\mathrm{lid}}(x_t^i)|^2
		+|u_\theta(x_b^i,0)|^2
		+|u_\theta(0,y_l^i)|^2
		+|u_\theta(1,y_r^i)|^2
		\Big],\\
		\mathcal{L}_{bc}^{v}(\theta)
		&=\frac{1}{N_{bc}}\sum_{i=1}^{N_{bc}}
		\Big[
		|v_\theta(x_t^i,1)|^2
		+|v_\theta(x_b^i,0)|^2
		+|v_\theta(0,y_l^i)|^2
		+|v_\theta(1,y_r^i)|^2
		\Big].
	\end{aligned}
\end{equation*}
Each training iteration uses a fresh set of $N_{bc}=100$ collocation points drawn from a uniform distribution along each wall. To alleviate the velocity discontinuities at the upper corners, the top-wall term in $\mathcal{L}_{bc}^{u}$ employs the smoothed lid-velocity profile
\begin{equation*}
	u_{\mathrm{lid}}(x)
	=
	1-\frac{\cosh\left(50(x-0.5)\right)}{\cosh(25)}.
\end{equation*}
The sampling settings are identical across all four methods. Grad Norm and NTK-based baselines are instead trained with the following loss:
\begin{equation*}
	\mathcal{L}(\theta)
	=\lambda_r\mathcal{L}_r(\theta)
	+\lambda_{bc}\mathcal{L}_{bc}(\theta).
\end{equation*}
Optimization spans $40{,}000$ iterations, with the learning rate decaying as $\eta_k=10^{-3}(0.9)^{k/1000}$. The proposed optimizer is implemented with $\delta_1=10^{-8}$. The relative $L^2$ errors are evaluated against the reference solution provided by Wang et al.~\cite{wang2023expertsguidetrainingphysicsinformed}:
\begin{equation*}
	\label{eq:ldc_relative_errors}
	e_q=
	\frac{\|q_\theta-q_{\mathrm{ref}}\|_2}
	{\|q_{\mathrm{ref}}\|_2},
	\qquad q\in\{u,v,U\},
\end{equation*}
where $U=\sqrt{u^2+v^2}$ is the velocity magnitude.

\begin{figure}[htbp]
	\centering
	\includegraphics[width=\textwidth]{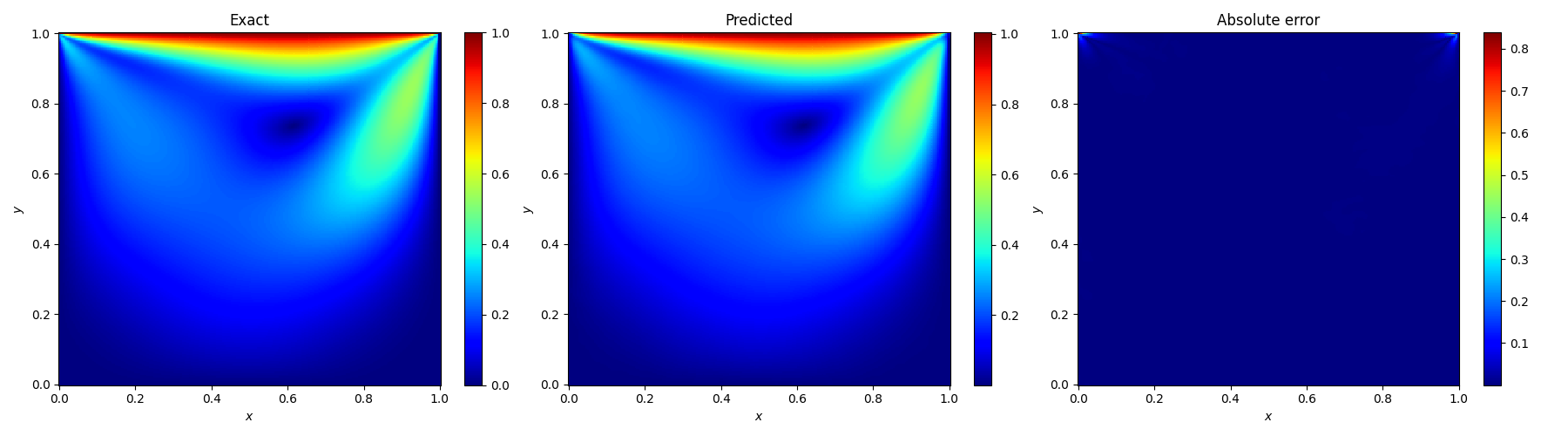}
	\caption{Lid-driven cavity flow at $Re=100$: reference velocity magnitude (left), prediction obtained with GMAdam (middle), and pointwise absolute error $|U_{\mathrm{ref}}-U_\theta|$ (right) after $40{,}000$ iterations.}
	\label{fig:ldc_re100_solution}
\end{figure}

As illustrated in Fig.~\ref{fig:ldc_re100_solution}, the predicted field reproduces the main features of the reference velocity-magnitude distribution, with relatively small pointwise errors across most of the cavity.

Table~\ref{tab:ldc_re100_comparison} indicates that the four methods attain comparable predictive accuracy. GMAdam improves upon AdaTask across all three error measures with similar training time, while achieving results close to those of Grad Norm. Although NTK-based weighting produces the lowest errors, its reduction in $e_U$ relative to the proposed method is approximately $1.7\%$, accompanied by an $85\%$ increase in training time.
\begin{table}[H]
	\centering
	\caption{Relative $L^2$ errors and training times for the lid-driven cavity problem at $Re=100$ after $40{,}000$ iterations.}
	\label{tab:ldc_re100_comparison}
	\begin{tabular}{lcccc}
		\toprule
		Method & $e_u$ & $e_v$ & $e_U$ & Training time (s) \\
		\midrule
		AdaTask
		& $4.3435\times10^{-2}$
		& $3.1404\times10^{-2}$
		& $3.7977\times10^{-2}$
		& 978.62 \\
		GMAdam
		& $4.3283\times10^{-2}$
		& $3.0456\times10^{-2}$
		& $3.7562\times10^{-2}$
		& 981.44 \\
		Grad-norm
		& $4.3502\times10^{-2}$
		& $2.8727\times10^{-2}$
		& $3.7503\times10^{-2}$
		& 970.59 \\
		NTK
		& $\mathbf{4.2995\times10^{-2}}$
		& $\mathbf{2.6276\times10^{-2}}$
		& $\mathbf{3.6914\times10^{-2}}$
		& 1816.74 \\
		\bottomrule
	\end{tabular}
\end{table}

\begin{figure}[H]
	\centering
	\begin{subfigure}[t]{0.24\textwidth}
		\centering
		\includegraphics[width=\linewidth]{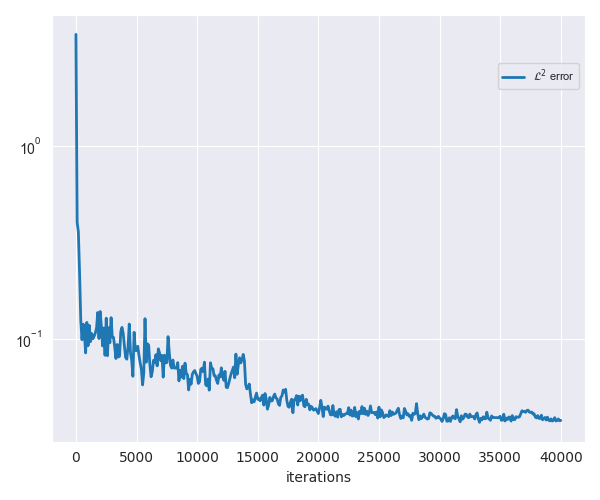}
		\caption{}
		\label{fig:ldc_re100_l2_adatask}
	\end{subfigure}
	\hfill
	\begin{subfigure}[t]{0.24\textwidth}
		\centering
		\includegraphics[width=\linewidth]{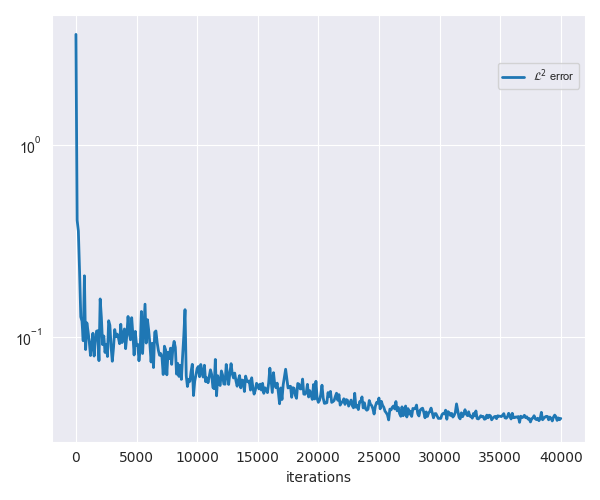}
		\caption{}
		\label{fig:ldc_re100_l2_ours}
	\end{subfigure}
	\hfill
	\begin{subfigure}[t]{0.24\textwidth}
		\centering
		\includegraphics[width=\linewidth]{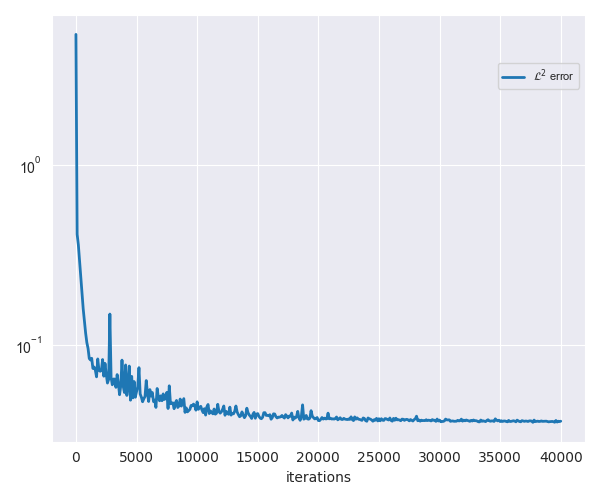}
		\caption{}
		\label{fig:ldc_re100_l2_gradnorm}
	\end{subfigure}
	\hfill
	\begin{subfigure}[t]{0.24\textwidth}
		\centering
		\includegraphics[width=\linewidth]{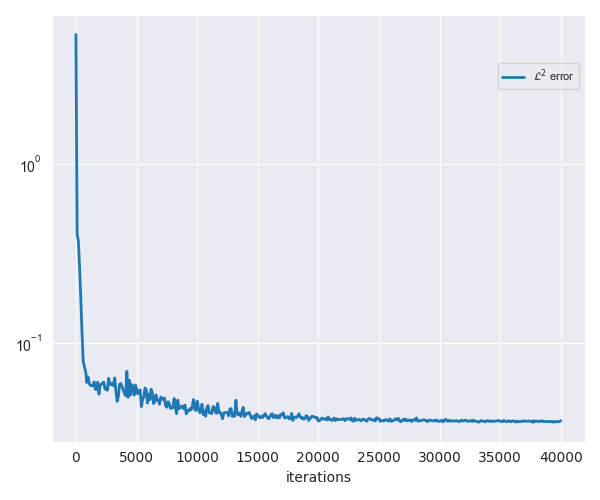}
		\caption{}
		\label{fig:ldc_re100_l2_ntk}
	\end{subfigure}
	\caption{Evolution of the relative $L^2$ error in velocity magnitude over $40{,}000$ training iterations at $Re=100$ for (a) AdaTask, (b) GMAdam, (c) Adam with Grad Norm, and (d) Adam with NTK-based loss weighting. The vertical axes are logarithmic.}
	\label{fig:ldc_re100_l2_histories}
\end{figure}

\begin{figure}[H]
	\centering
	\includegraphics[width=\textwidth]{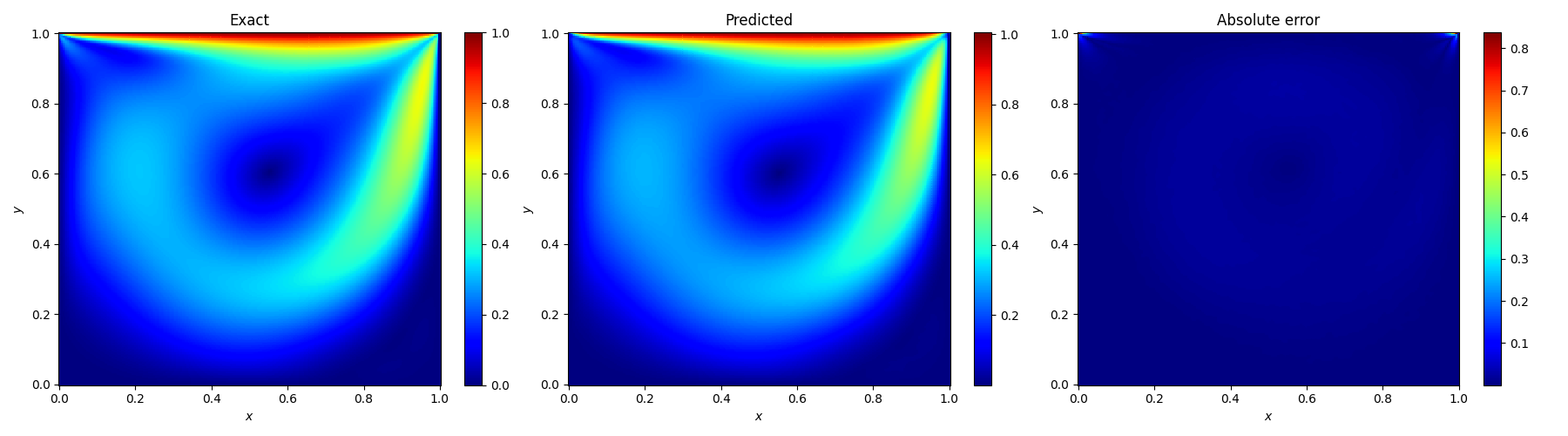}
	\caption{Lid-driven cavity flow at $Re=400$: reference velocity magnitude (left), prediction obtained with GMAdam (middle), and pointwise absolute error $|U_{\mathrm{ref}}-U_\theta|$ (right) after $80{,}000$ iterations.}
	\label{fig:ldc_re400_solution}
\end{figure}

As shown in Fig.~\ref{fig:ldc_re100_l2_histories}, the four methods exhibit different convergence behaviors despite their similar final accuracy. Gradient-norm-based and NTK-based weighting achieve relatively rapid initial error reduction, whereas AdaTask and the proposed method experience more pronounced transient fluctuations. The subsequent flattening of the error trajectories indicates diminishing improvements during the later stages of training.
The $Re=400$ experiment follows the same numerical setup as above, with the Fourier embedding scale adjusted to $\sigma=4$. To better assess the late-stage convergence behavior of the four methods, the training horizon is extended to $80{,}000$ iterations. As the Reynolds number increases to $400$, the interior low-velocity region shifts toward the cavity center, a feature that is well captured by the proposed method, as shown in Fig.~\ref{fig:ldc_re400_solution}.

Table~\ref{tab:ldc_re400_comparison} reveals a substantially greater advantage of the proposed method over AdaTask at $Re=400$ than at $Re=100$. In particular, the proposed method reduces $e_U$ by approximately $77\%$ relative to AdaTask, accompanied by marked reductions in both velocity-component errors, while leaving the training time virtually unchanged. Although Grad Norm and NTK-based weighting attain lower errors, only the former does so at a comparable computational cost.
\begin{table}[H]
	\centering
	\caption{Relative $L^2$ errors and training times for the lid-driven cavity problem at $Re=400$ after $80{,}000$ iterations.}
	\label{tab:ldc_re400_comparison}
	\begin{tabular}{lcccc}
		\toprule
		Method & $e_u$ & $e_v$ & $e_U$ & Training time (s) \\
		\midrule
		AdaTask
		& $2.6606\times10^{-1}$
		& $2.6578\times10^{-1}$
		& $2.3981\times10^{-1}$
		& 1893.92 \\
		GMAdam
		& $5.9259\times10^{-2}$
		& $5.2851\times10^{-2}$
		& $5.5362\times10^{-2}$
		& 1899.63 \\
		Grad Norm
		& $\mathbf{4.2040\times10^{-2}}$
		& $\mathbf{1.9143\times10^{-2}}$
		& $\mathbf{3.3066\times10^{-2}}$
		& 1862.50 \\
		NTK
		& $4.4458\times10^{-2}$
		& $2.5805\times10^{-2}$
		& $3.6427\times10^{-2}$
		& 3480.67 \\
		\bottomrule
	\end{tabular}
\end{table}

\begin{figure}[H]
	\centering
	\begin{subfigure}[t]{0.24\textwidth}
		\centering
		\includegraphics[width=\linewidth]{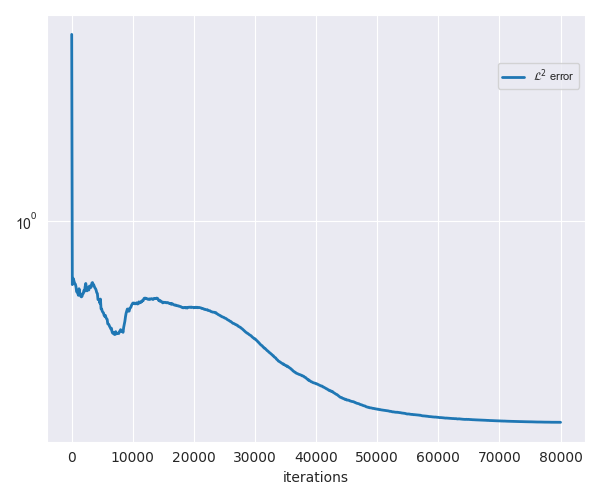}
		\caption{}
		\label{fig:ldc_re400_l2_adatask}
	\end{subfigure}
	\hfill
	\begin{subfigure}[t]{0.24\textwidth}
		\centering
		\includegraphics[width=\linewidth]{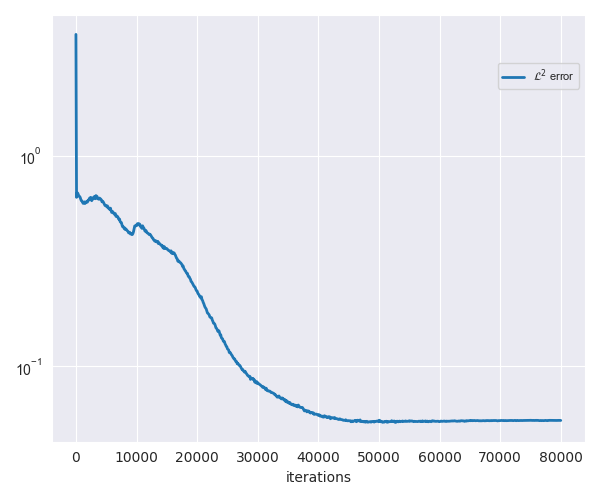}
		\caption{}
		\label{fig:ldc_re400_l2_ours}
	\end{subfigure}
	\hfill
	\begin{subfigure}[t]{0.24\textwidth}
		\centering
		\includegraphics[width=\linewidth]{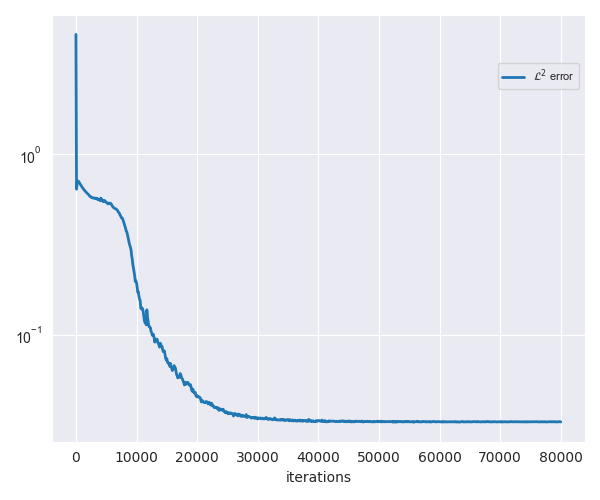}
		\caption{}
		\label{fig:ldc_re400_l2_gradnorm}
	\end{subfigure}
	\hfill
	\begin{subfigure}[t]{0.24\textwidth}
		\centering
		\includegraphics[width=\linewidth]{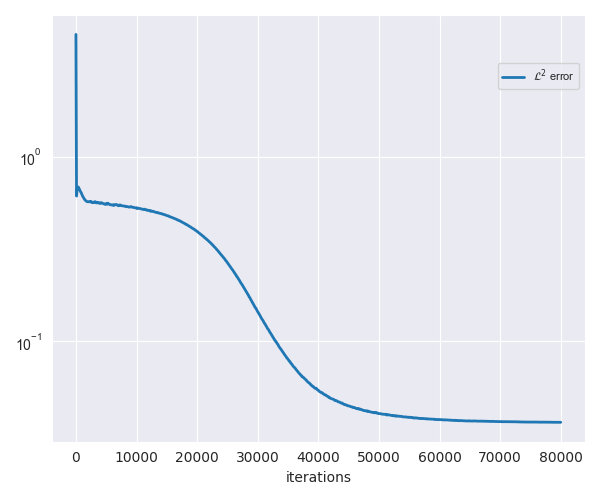}
		\caption{}
		\label{fig:ldc_re400_l2_ntk}
	\end{subfigure}
	\caption{Evolution of the relative $L^2$ error in velocity magnitude over $80{,}000$ training iterations at $Re=400$ for (a) AdaTask, (b) GMAdam, (c) Adam with Grad Norm, and (d) Adam with NTK-based loss weighting. The vertical axes are logarithmic.}
	\label{fig:ldc_re400_l2_histories}
\end{figure}

Beyond the final-error comparison, the convergence histories in Fig.~\ref{fig:ldc_re400_l2_histories} reveal that, following transient increases in the relative $L^2$ error, GMAdam returns to a sustained downward trend earlier than AdaTask and achieves a more pronounced subsequent reduction.
\bibliographystyle{unsrt}  
\bibliography{references}

\begin{thebibliography}{10}

\bibitem{E_2017}
Weinan E, Jiequn Han, and Arnulf Jentzen.
\newblock Deep learning-based numerical methods for high-dimensional parabolic
  partial differential equations and backward stochastic differential
  equations.
\newblock {\em Communications in Mathematics and Statistics}, 5(4):349--380,
  November 2017.

\bibitem{han2018solving}
Jiequn Han, Arnulf Jentzen, and Weinan E.
\newblock Solving high-dimensional partial differential equations using deep
  learning.
\newblock {\em Proceedings of the National Academy of Sciences},
  115(34):8505--8510, 2018.

\bibitem{e2017deepritzmethoddeep}
Weinan E and Bing Yu.
\newblock The deep ritz method: A deep learning-based numerical algorithm for
  solving variational problems, 2017.

\bibitem{liao2019deep}
Yulei Liao and Pingbing Ming.
\newblock Deep nitsche method: Deep ritz method with essential boundary
  conditions.
\newblock {\em arXiv preprint arXiv:1912.01309}, 2019.

\bibitem{zang2020weak}
Yaohua Zang, Gang Bao, Xiaojing Ye, and Haomin Zhou.
\newblock Weak adversarial networks for high-dimensional partial differential
  equations.
\newblock {\em Journal of Computational Physics}, 411:109409, 2020.

\bibitem{raissi2019physics}
Maziar Raissi, Paris Perdikaris, and George~E Karniadakis.
\newblock Physics-informed neural networks: A deep learning framework for
  solving forward and inverse problems involving nonlinear partial differential
  equations.
\newblock {\em Journal of Computational physics}, 378:686--707, 2019.

\bibitem{sun2020surrogate}
Luning Sun, Han Gao, Shaowu Pan, and Jian-Xun Wang.
\newblock Surrogate modeling for fluid flows based on physics-constrained deep
  learning without simulation data.
\newblock {\em Computer Methods in Applied Mechanics and Engineering},
  361:112732, 2020.

\bibitem{cai2021physics}
Shengze Cai, Zhiping Mao, Zhicheng Wang, Minglang Yin, and George~Em
  Karniadakis.
\newblock Physics-informed neural networks (pinns) for fluid mechanics: A
  review.
\newblock {\em Acta Mechanica Sinica}, 37(12):1727--1738, 2021.

\bibitem{eivazi2022physics}
Hamidreza Eivazi, Mojtaba Tahani, Philipp Schlatter, and Ricardo Vinuesa.
\newblock Physics-informed neural networks for solving reynolds-averaged
  navier--stokes equations.
\newblock {\em Physics of Fluids}, 34(7), 2022.

\bibitem{haghighat2021sciann}
Ehsan Haghighat and Ruben Juanes.
\newblock Sciann: A keras/tensorflow wrapper for scientific computations and
  physics-informed deep learning using artificial neural networks.
\newblock {\em Computer Methods in Applied Mechanics and Engineering},
  373:113552, 2021.

\bibitem{smith2022hyposvi}
Jonthan~D Smith, Zachary~E Ross, Kamyar Azizzadenesheli, and Jack~B Muir.
\newblock Hyposvi: Hypocentre inversion with stein variational inference and
  physics informed neural networks.
\newblock {\em Geophysical Journal International}, 228(1):698--710, 2022.

\bibitem{kissas2020machine}
Georgios Kissas, Yibo Yang, Eileen Hwuang, Walter~R Witschey, John~A Detre, and
  Paris Perdikaris.
\newblock Machine learning in cardiovascular flows modeling: Predicting
  arterial blood pressure from non-invasive 4d flow mri data using
  physics-informed neural networks.
\newblock {\em Computer methods in applied mechanics and engineering},
  358:112623, 2020.

\bibitem{sahli2020physics}
Francisco Sahli~Costabal, Yibo Yang, Paris Perdikaris, Daniel~E Hurtado, and
  Ellen Kuhl.
\newblock Physics-informed neural networks for cardiac activation mapping.
\newblock {\em Frontiers in Physics}, 8:42, 2020.

\bibitem{9403414}
Zhiwei Fang.
\newblock A high-efficient hybrid physics-informed neural networks based on
  convolutional neural network.
\newblock {\em IEEE Transactions on Neural Networks and Learning Systems},
  33(10):5514--5526, 2022.

\bibitem{kovacs2022conditional}
Alexander Kovacs, Lukas Exl, Alexander Kornell, Johann Fischbacher, Markus
  Hovorka, Markus Gusenbauer, Leoni Breth, Harald Oezelt, Masao Yano, Noritsugu
  Sakuma, et~al.
\newblock Conditional physics informed neural networks.
\newblock {\em Communications in Nonlinear Science and Numerical Simulation},
  104:106041, 2022.

\bibitem{liu2019multi}
Dehao Liu and Yan Wang.
\newblock Multi-fidelity physics-constrained neural network and its application
  in materials modeling.
\newblock {\em Journal of Mechanical Design}, 141(12):121403, 2019.

\bibitem{chen2020physics}
Yuyao Chen, Lu~Lu, George~Em Karniadakis, and Luca Dal~Negro.
\newblock Physics-informed neural networks for inverse problems in nano-optics
  and metamaterials.
\newblock {\em Optics express}, 28(8):11618--11633, 2020.

\bibitem{wang2025gradient}
Sifan Wang, Ananyae~Kumar Bhartari, Bowen Li, and Paris Perdikaris.
\newblock Gradient alignment in physics-informed neural networks: A
  second-order optimization perspective.
\newblock {\em arXiv preprint arXiv:2502.00604}, 2025.

\bibitem{wang2021understanding}
Sifan Wang, Yujun Teng, and Paris Perdikaris.
\newblock Understanding and mitigating gradient flow pathologies in
  physics-informed neural networks.
\newblock {\em SIAM Journal on Scientific Computing}, 43(5):A3055--A3081, 2021.

\bibitem{wang2022and}
Sifan Wang, Xinling Yu, and Paris Perdikaris.
\newblock When and why pinns fail to train: A neural tangent kernel
  perspective.
\newblock {\em Journal of Computational Physics}, 449:110768, 2022.

\bibitem{jacot2018neural}
Arthur Jacot, Franck Gabriel, and Cl{\'e}ment Hongler.
\newblock Neural tangent kernel: Convergence and generalization in neural
  networks.
\newblock {\em Advances in neural information processing systems}, 31, 2018.

\bibitem{AdaTask_AAAI2023}
Enneng Yang, Junwei Pan, Ximei Wang, Haibin Yu, Li~Shen, Xihua Chen, Lei Xiao,
  Jie Jiang, and Guibing Guo.
\newblock Adatask: A task-aware adaptive learning rate approach to multi-task
  learning.
\newblock {\em Proceedings of the AAAI Conference on Artificial Intelligence},
  37(9):10745--10753, 2023.

\bibitem{pmlr-v202-yao23c}
Jiachen Yao, Chang Su, Zhongkai Hao, Songming Liu, Hang Su, and Jun Zhu.
\newblock {M}ulti{A}dam: Parameter-wise scale-invariant optimizer for
  multiscale training of physics-informed neural networks.
\newblock In Andreas Krause, Emma Brunskill, Kyunghyun Cho, Barbara Engelhardt,
  Sivan Sabato, and Jonathan Scarlett, editors, {\em Proceedings of the 40th
  International Conference on Machine Learning}, volume 202 of {\em Proceedings
  of Machine Learning Research}, pages 39702--39721. PMLR, 23--29 Jul 2023.

\bibitem{barakat2021convergence}
Anas Barakat and Pascal Bianchi.
\newblock Convergence and dynamical behavior of the adam algorithm for
  nonconvex stochastic optimization.
\newblock {\em SIAM Journal on Optimization}, 31(1):244--274, 2021.

\bibitem{pmlr-v145-ma22a}
Chao Ma, Lei Wu, and Weinan E.
\newblock A qualitative study of the dynamic behavior for adaptive gradient
  algorithms.
\newblock In Joan Bruna, Jan Hesthaven, and Lenka Zdeborova, editors, {\em
  Proceedings of the 2nd Mathematical and Scientific Machine Learning
  Conference}, volume 145 of {\em Proceedings of Machine Learning Research},
  pages 671--692. PMLR, 16--19 Aug 2022.

\bibitem{glorot2010understanding}
Xavier Glorot and Yoshua Bengio.
\newblock Understanding the difficulty of training deep feedforward neural
  networks.
\newblock In {\em Proceedings of the thirteenth international conference on
  artificial intelligence and statistics}, pages 249--256. JMLR Workshop and
  Conference Proceedings, 2010.

\bibitem{kingma2015iclr-adam}
Diederik~P. Kingma and Jimmy Ba.
\newblock {Adam: A Method for Stochastic Optimization}.
\newblock In {\em International Conference on Learning Representations}, 2015.

\bibitem{wang2023expertsguidetrainingphysicsinformed}
Sifan Wang, Shyam Sankaran, Hanwen Wang, and Paris Perdikaris.
\newblock An expert's guide to training physics-informed neural networks, 2023.

\bibitem{wang2021eigenvector}
Sifan Wang, Hanwen Wang, and Paris Perdikaris.
\newblock On the eigenvector bias of fourier feature networks: From regression
  to solving multi-scale pdes with physics-informed neural networks.
\newblock {\em Computer Methods in Applied Mechanics and Engineering},
  384:113938, 2021.

\bibitem{pmlr-v306-bai26o}
Zhiwei Bai, Zhangchen Zhou, Jiajie Zhao, Xiaolong Li, Zhiyu Li, Feiyu Xiong,
  Hongkang Yang, Yaoyu Zhang, and Zhi-Qin~John Xu.
\newblock Adaptive preconditioners trigger loss spikes in {A}dam.
\newblock In Tong Zhang, Miroslav Dudik, Martin Jaggi, Alekh Agarwal, Sharon
  Li, Dale Schuurmans, Jerry Zhu, Felix Berkenkamp, Hanze Dong, and Alberto
  Bietti, editors, {\em Proceedings of the 43rd International Conference on
  Machine Learning}, volume 306 of {\em Proceedings of Machine Learning
  Research}, pages 5566--5601. PMLR, 06--11 Jul 2026.

\bibitem{zhang2026adamconvergesmodificationupdate}
Yushun Zhang, Bingran Li, Congliang Chen, Zhi-Quan Luo, and Ruoyu Sun.
\newblock Adam converges without any modification on update rules, 2026.

\bibitem{FENICHEL197953}
Neil Fenichel.
\newblock Geometric singular perturbation theory for ordinary differential
  equations.
\newblock {\em Journal of Differential Equations}, 31(1):53--98, 1979.

\bibitem{kuehn2015multiple}
Christian Kuehn et~al.
\newblock {\em Multiple time scale dynamics}, volume 191.
\newblock Springer, 2015.

\bibitem{Stein01051987}
Michael Stein.
\newblock Large sample properties of simulations using latin hypercube
  sampling.
\newblock {\em Technometrics}, 29(2):143--151, 1987.

\bibitem{10.1007/978-3-032-06774-6_1}
Nabil Jabareen, Dongsheng Yuan, Dingming Liu, Foo-Wei Ten, and S{\"o}ren
  Lukassen.
\newblock Anisotropic {Fourier} features for positional encoding in medical
  imaging.
\newblock In Christian Wachinger, Gijs Luijten, Shireen Elhabian, Karthik
  Gopinath, and Jan Egger, editors, {\em Shape in Medical Imaging}, pages
  1--13, Cham, 2026. Springer Nature Switzerland.

\end{thebibliography}

\end{document}